\documentclass[letterpaper, journal]{IEEEtran}
\usepackage{epsfig} % for 
\usepackage{varwidth}
\usepackage[abs]{overpic}
\usepackage{algorithmic}

\newcommand{\squeezeup}{\vspace{-2.5mm}}

\usepackage{cite}
\usepackage{enumerate}

\usepackage{epstopdf}

\ifCLASSINFOpdf
\else
\fi
\usepackage{amsmath}
\usepackage{amssymb}  % assumes amsmath package installed

\usepackage{algorithmic}
\usepackage{caption}
\usepackage[caption=false,font=footnotesize]{subfig}
\begin{document}

\IEEEpubid{Submitted on Sep-01-2016 to IEEE Transactions on Control Systems Technology. Manuscript ID: TCST-2016-0820.
%\\ Revised Manuscript ID: TCST-2017-0471. Rejected on Nov-10-2017.
}% Always comment-out the \thispagestyle{empty} command
\IEEEpubidadjcol

\title{Geometric Surface-Based Tracking Control of a Quadrotor UAV under Actuator Constraints
% for aggressive maneuvers
 }

\author{Michalis~Ramp$^1$~and~Evangelos~Papadopoulos$^2$% <-this % stops a space
\thanks{$^1$M. Ramp is with the Department of Mechanical Engineering, National Technical University of Athens, (NTUA) 15780 Athens, Greece.
{\tt\small rampmich@mail.ntua.gr}}% <-this % stops a space
\thanks{$^2$E. Papadopoulos is with the Department of Mechanical Engineering, NTUA, 15780 Athens (tel: +30-210-772-1440; fax: +30-210-772-1455).
{\tt\small egpapado@central.ntua.gr}}% <-this % stops a space
%\thanks{Manuscript received January 19, 2016; revised March 26, 2016.}
}

\maketitle

\begin{abstract}
This paper presents contributions on nonlinear tracking control systems for a quadrotor unmanned micro aerial vehicle.
New controllers are proposed based on nonlinear surfaces composed by tracking errors that evolve directly on the nonlinear configuration manifold thus inherently including in the control design the nonlinear characteristics of the SE(3) configuration space.
In particular geometric surface-based controllers are developed, and through rigorous stability proofs they are shown to have desirable closed loop properties that are almost global.
A region of attraction, independent of the position error, is produced and its effects are analyzed.
% Furthermore the controllers are supported by rigorous stability proofs.
%The final contribution that differentiates this work from existing geometric quadrotor controllers is the introduction of a strategy for the inclusion of a lower bound constraint on the controller generated thrusts, without interfering with the control objective, ensuring that the generated thrusts stay positive.
A strategy allowing the quadrotor to achieve precise attitude tracking while simultaneously following a desired position command and complying to actuator constraints in a computationally inexpensive manner is derived.
This important contribution differentiates this work from existing Geometric Nonlinear Control System solutions (GNCSs) since the commanded thrusts can be realized by the majority of quadrotors produced by the industry.
%In contrast, existing Geometric Nonlinear Control System solutions (GNCSs) produce negative thrusts during aggressive maneuvers that cannot be generated on a standard quadrotor.
The new features of the proposed GNCSs are illustrated by numerical simulations of aggressive maneuvers and a comparison with a GNCSs from the bibliography.
\end{abstract}


% \fi
\IEEEpeerreviewmaketitle

\section{Introduction}
\IEEEPARstart{Q}{uadrotor} unmanned aerial vehicles are characterized by a simple mechanical structure comprised of two pairs of counter rotating outrunner motors where each one is driving a dedicated propeller, resulting to a platform with high thrust-to-weight ratio, able to achieve vertical takeoff and landing (VTOL) maneuvers and operate in a broad spectrum of flight scenarios.
Furthermore the quadrotor as a platform has been identified to have good flight endurance characteristics and acceptable payload transporting potential for a plethora of applications \cite{Loianno}.
Although the quadrotor UAV has six degrees of freedom, it is underactuated since it has only four inputs.
As a result the quadrotor can only track four commands or less.
Even though the research of quadrotor unmanned aerial vehicles is in a mature state, the platform's current popularity is at an all time high with the vehicle to be easily accessible to the general public in a large variety of configurations and features from a number of vendors.

A plethora of control systems, linear and nonlinear, have been proposed for this platform.
In \cite{Castillo} a linear controller was developed based on Lyapunov analysis, and experimental results were presented demonstrating the vehicle performing vertical takeoff, steady hover and landing, autonomously.

%To address issues that arise when deviating significantly from the typical hover flight regime, three aerodynamic effects were investigated and the results were used to improve the performance of a type of linear-PID controller that was employed \cite{Hoffmann07}.
%The results however revealed that existing models and control techniques were inadequate for accurate trajectory tracking at high speeds and in uncontrolled environments.

A back-flip aerobatic maneuver and a decentralized collision avoidance algorithm for multiple quadrotors was designed using hybrid decomposition and reachable set theory \cite{Gillula}.
However this type of analysis only applies to a limited flight envelope.
%In [] a linear controller was were developed to augment the stability properties of the platform usually in near hover conditions.
The controller design and the trajectory generation for a quadrotor maneuvering in a tightly constrained environment was addressed using differential flatness, \cite{Mellinger}.
This resulted in the development of an algorithm that enables real time generation of optimal trajectories by the minimization of cost functionals derived from the norm of the fourth derivative of the position.
The algorithm enables safe passage through corridors while satisfying constraints on velocities, accelerations, and inputs, although no stability proof was given.

%An integral predictive and robust nonlinear H$_{\infty}$ control strategy addressed the path following problem for a quadrotor helicopter, \cite{Raffo}.
%The control strategy took under consideration the presence of aerodynamic disturbances and the parametric and structural uncertainties acting on all degrees of freedom.
%A simulation model taking into account the variation of the aerodynamical coefficients due to vehicle motion was introduced, and a controller using Integral Backstepping was proposed for full control of the quadrotor attitude, altitude and position \cite{Siegwart}. Results showcased autonomous take-off, hover, landing and collision avoidance.

A fuzzy controller was applied to control the position and orientation of a quadrotor with acceptable performance in simulation \cite{Matilde}.
However, due to the interdependence of variables, the tuning of the controller parameters presented a considerable challenge and was done through trial and error.

A nonlinear dynamic model in a form suited for backstepping control was presented in \cite{Madani}, followed by the design of a backstepping control law based on Lyapunov stability theory by means of decomposition of the system equations into three interconnected subsystems.
However the simulation results involved only tracking of a basic step trajectory of a smooth transition to a desired position and yaw angle.

In most of the aforementioned works, and in the majority of those in the literature, the design of controllers is based on minimum attitude representations such as the Euler angles, i.e. on a local representation of attitude.
Thus, they involve convoluted, long and complicated expressions, they exhibit singularities during large angle rotational maneuvers, and they restrict significantly the quadrotors operational envelope.

\IEEEpubidadjcol

A controller based on quaternions was developed in \cite{Tayebi}.
Although quaternions are singularity-free, the three-sphere S$^3$ quaternion space double-covers SO(3), meaning that a quaternion control system entails convergence to either of the two disconnected antipodal points of S$^3$, both representing the same global attitude \cite{Mayhew}.
As a result, if this ambiguity is not dealt with during control design, the quadrotor becomes sensitive to noise, the quaternion controller can become discontinuous and can give rise to unwinding phenomena where the rigid body unnecessarily rotates even though its attitude is extremely close to the desired orientation \cite{MayhewNon},\cite{Chaturvedi}.

The quadrotor as a dynamic system evolves on a nonlinear manifold;
Hence, it cannot be described globally with Euclidean spaces.
Recently, attitude control was studied within a geometric framework where the dynamics of a rigid body was globally expressed in SO(3) or S$^2$.
Control systems developed in a geometric framework encompass the attributes of the system nonlinear manifolds in the characterization of the configuration manifold, an additional advantage, avoiding ambiguities and singularities associated with minimal representations of attitude.
%By using geometric control techniques, control systems are developed that inherently encompass the attributes of the systems nonlinear manifolds in the characterization of the configuration manifold, avoiding in this manner ambiguities and singularities associated with minimal representations of attitude.
This approach was applied to fully actuated and under-actuated dynamic systems on Lie groups to achieve almost global asymptotic stability \cite{Chaturvedi} -\cite{obstruction}.

Towards this direction, the dynamics of a quadrotor UAV was globally expressed on the special Euclidean group SE(3) and nonlinear controllers, designed directly on the nonlinear configuration manifold, were developed with flight maneuvers defined by a concatenation of three flight modes, an attitude, a velocity and a position flight mode \cite{geomquadlee}.
The resulting Geometric Nonlinear Control System solutions (GNCSs) were shown to have desirable properties that are almost global in each mode illustrating the versatility and generality of the proposed approach/solution.
The aforementioned results were extended by a robust GNCS proposed in \cite{geomquadlee_asian}, by adaptive GNCSs in \cite{qeoadaplee},\cite{qeoadapfar}, and by a nonlinear PID GNCS in \cite{qeopidfar}.
%However in \cite{geomquadlee},\cite{geomquadlee_asian},\cite{qeopidfar}, the GNCS-generated thrusts become negative during aggressive maneuvers which can be produced only by a nonstandard quadrotor equipped with variable pitch propellers (or other additional means).
However the GNCS-generated thrusts become negative during aggressive maneuvers in \cite{geomquadlee},\cite{geomquadlee_asian},\cite{qeopidfar}.
Consequently because negative thrusts can be produced only by a nonstandard quadrotor equipped with variable pitch propellers (or other additional means), the performance assurances produced in \cite{geomquadlee},\cite{geomquadlee_asian},\cite{qeopidfar} are valid for a standard quadrotor only when the controller generated thrusts are realizable by a typical outrunner motor.

In the same geometric context, the quadrotor as means of load transportation was investigated.
The dynamics and control of quadrotor(s) with a payload that is connected via flexible cable(s) was studied in \cite{qeocablefar} -\cite{qeocablesfar}.
The problem of the stabilization of a rigid body payload with multiple cooperative quadrotors was addressed, with the theoretical results to be supported by an experimental implementation \cite{stabfar}.

This paper follows the geometric framework.
A GNCS for a quadrotor UAV is developed directly on the special Euclidean group thus inherently entailing in the control design the characteristics of the nonlinear configuration manifold and avoiding singularities and ambiguities associated with minimal attitude representations.
The key contributions of this work are: 

(1) We propose controllers (an attitude and a position controller) based on nonlinear surfaces composed by tracking errors that evolve directly on the nonlinear configuration manifold.
These controllers allow for precision pose tracking by tuning three gains per controller and are able to follow an attitude tracking command and a position tracking command.
%Moreover it is showed in simulation that our proposed controller can operate in lesser control loop frequencies in comparison to a benchmark geometric controller in the bibliography.

(2) In contrast to \cite{geomquadlee} -\cite{qeopidfar}, we produce rigorous stability proofs and regions of attraction both with and without restrictions on the initial position/velocity error.
In both cases it is shown that the position controller structure is characterized by almost global exponential attractiveness.
%Furthermore we compare the three regions of exponential stability and investigate the advantages that they offer in regards to trajectory design.
It is shown that the basin of attraction that does not depend explicitly on the initial position/velocity error is smaller with respect to other basins, but in contrast to these it is a function of only \textit{the control gains} and \textit{the quadrotor mass}, introducing simplicity in trajectory design.

(3) A strategy allowing the quadrotor to achieve, (i) precise attitude tracking, while (ii) following a desired position command, and (iii) complying to actuator constraints in a computationally inexpensive manner is developed.
%This is an unprecedented contribution differentiating this work from existing Geometric Nonlinear Control System solutions (GNCSs) since the generated thrusts can be realized by the majority of quadrotors produced by the industry.
%In contrast, existing Geometric Nonlinear Control System solutions (GNCSs) produce negative thrusts during aggressive maneuvers that cannot be generated on a standard quadrotor.The final contribution that differentiates this work from existing GNCSs is a strategy for the inclusion of a lower bound constraint on the controller generated thrusts, without interfering with the control objective, ensuring that the generated thrusts stay positive.
As a result, enhanced capabilities are achieved, while the generated thrusts are realizable by the majority of existing quadrotors.
%Furthermore the performance assurances that are extracted by the rigorous mathematical proofs can be extended to an actual quadrotor since through the strategy, we account for this important factor.
Thus the performance assurances that are extracted by the rigorous mathematical proofs can be extended to an actual quadrotor.

This strategy in conjunction with the developed controllers in (1) and the regions of attraction from (2) comprise our GNCS.
% for a quadrotor UAV.
%This comes in contrast with existing geometric nonlinear tracking control systems where during aggressive maneuvers they require negative thrusts that can only be produced by a nonstandard quadrotor equipped with variable pitch propellers (or other means), see the simulation results for the controller generated thrusts in \cite{geomquadlee},\cite{geomquadlee_asian},\cite{qeopidfar}.
The proposed strategies are validated in simulation.
To the authors best knowledge, 
%the design of geometric surface-based control systems for a quadrotor UAV, the production of position free region of attraction and the introduction of a strategy for precise attitude tracking, while simultaneously following a desired position command, and comply to actuator constraints in a computationally cheap manner for a quadrotor UAV on SE(3) is unprecedented.
the above contributions are completely novel and extend the quadrotor UAV nonlinear control methodologies on SE(3).
\section{Quadrotor Kinetics Model}
\begin{figure}[thpb]
%\begin{figure}[!h]
      \centering
      \includegraphics[width=1\columnwidth]{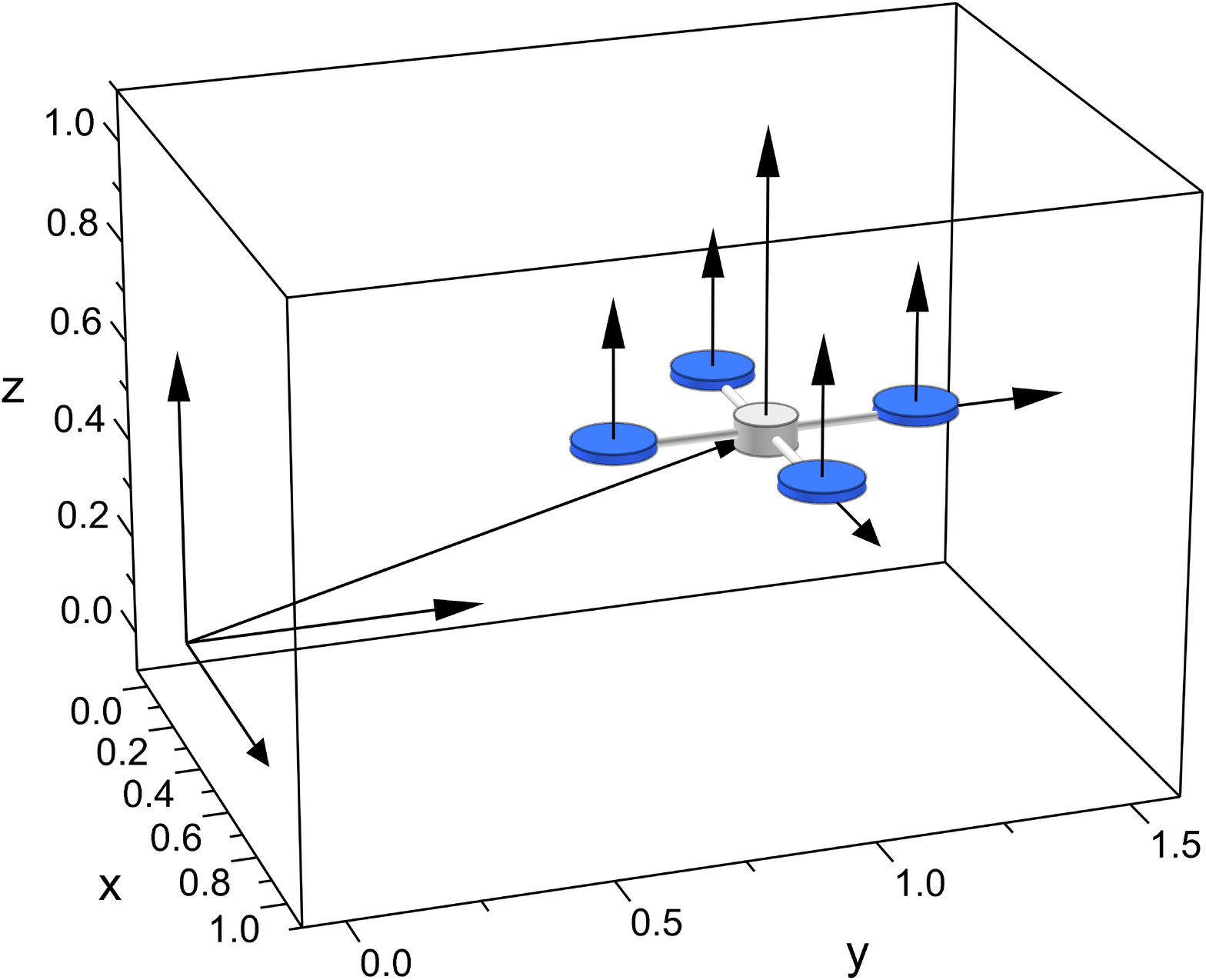}
      % Body-fixed
      \put(-82,74){\parbox{\columnwidth}{$\mathbf{e}_{1}$}}
      \put(-59,101){\parbox{0.4\linewidth}{$\mathbf{e}_{2}$}}
      \put(-101,142){\parbox{0.4\linewidth}{$\mathbf{e}_{3}$}}
      % Position Vector
      \put(-143,71){\parbox{0.4\linewidth}{$\mathbf{x}$}}
      % Inertial
      \put(-177,37){\parbox{\columnwidth}{$\mathbf{E}_{1}$}}
      %\put(-145,64){\parbox{0.4\linewidth}{$\mathbf{E}_{2}$}}
      \put(-147,58){\parbox{0.4\linewidth}{$\mathbf{E}_{2}$}}
      \put(-196,107){\parbox{0.4\linewidth}{$\mathbf{E}_{3}$}}
      % Thrusts
      \put(-95,112){\parbox{0.4\linewidth}{$f_{1}$}}
      \put(-82,122){\parbox{0.4\linewidth}{$f_{2}$}}
      \put(-113,127){\parbox{0.4\linewidth}{$f_{3}$}}
      \put(-129,116){\parbox{0.4\linewidth}{$f_{4}$}}
      \caption{Quadrotor with the coordinate frames, unit vectors that define it and actuator forces.}
\label{Quadrotor}
%\end{picture}
\end{figure}

The quadrotor studied is comprised by two pairs of counter rotating out-runner motors see Fig. \ref{Quadrotor}.
Each motor drives a dedicated propeller and generates thrust and torque normal to the plane produced by the centers of mass (CM) of the four rotors. 
An inertial reference frame  I$_{R}\big\{\mathbf{E}_1,\mathbf{E}_2,\mathbf{E}_3\big\}$ and a body-fixed frame I$_{b}\big\{\mathbf{e}_1,\mathbf{e}_2,\mathbf{e}_3\big\}$ are employed with the origin of the latter to be located at the quadrotor CM. 
The first two axes of I$_{b}$ are co-linear with the two quadrotor legs as depicted in Fig. \ref{Quadrotor} and lie on the same plane defined by the CM of the four rotors and the CM of the quadrotor.
%All main quadrotor properties, and all kinematic quantities are displayed in Table \ref{Table}. 
%\begin{overpic}
%[width=1\columnwidth,grid,scale=0.5,unit=0.5mm]{./figures/Quad.eps}
%%[width=1\columnwidth]{./figures/Quad.eps}
%\put(110,78){\parbox{0.4\linewidth}{$f_{1}$}}
%\end{overpic}

The following apply throughout the paper.
The thrust of each propeller is considered to be the actual control input and acts along the direction of the propeller axis which is co-linear with the $\mathbf{e}_{3}$ body-fixed axis.
The first and third propellers generate positive thrust along the direction of $\mathbf{e}_{3}$ when rotating clockwise, while the second and fourth propellers generate positive thrust along the direction of $\mathbf{e}_{3}$ when rotating counterclockwise.
The magnitude of the total thrust is denoted by $f=\sum_{i=1}^{4} f_i\in\mathbb{R}$ and it is positive when it acts along $\mathbf{e}_{3}$ and negative when acts along the 
$-\mathbf{e}_{3}$ direction, where $f_{i}$ and other system variables are defined in Table \ref{Table}.

\begin{table}[h]
%\captionsetup{justification=centering}
\caption{Definitions of variables}
\label{Table}
\begin{center}
\begin{tabular}{|l|l|}
\hline
$\mathbf{x}\in\mathbb{R}^3$ & Quadrotor CM position wrt. I$_R$ in I$_R$\\
\hline
$\mathbf{v}\in\mathbb{R}^3$ & Quadrotor CM velocity wrt. I$_R$ in I$_R$\\
\hline
$^{b}\boldsymbol{\omega}\in\mathbb{R}^3$ & Quadrotor angular velocity wrt I$_R$ in I$_{b}$\\
\hline
$\mathbf{R}\in\text{SO}\left(3\right)$ & Rotation matrix from $\mathbf{I}_b$ to $\mathbf{I}_R$ frame\\
\hline
$^b\mathbf{u}\in\mathbb{R}^3$ & Controller generated torque $^b\mathbf{u}{=}[{}^{b}u_{1};{}^{b}u_{2};{}^{b}u_{3}]$ in I$_b$\\
\hline
$f_i\in\mathbb{R}$ & Force produced by the i-th propeller along $e_{3}$ \\
\hline
$b_T\in\mathbb{R}^{+}$ & Torque coefficient \\%where $(-1)^{i}b_{T}f_{i}e_{3}$ the\\& torque produced by the i-th propeller\\
\hline
$g\in\mathbb{R}$ & Gravity constant\\
\hline
$d\in\mathbb{R}^{+}$ & Distance between system CM and each motor axis\\
\hline
$\mathbf{J}\in\mathbb{R}^{3\times3}$ & Inertial matrix (IM) of the quadrotor in I$_b$\\
\hline
$m\in\mathbb{R}$ & Quadrotor total mass\\
\hline
$\lambda_{min,max}(.)$ & Minimum, maximum eigenvalue of $(.)$ respectively\\
\hline
\end{tabular}
\end{center}
\end{table}

The motor torques corresponding to each propeller are assumed to be proportional to each propeller thrust.
Specifically, the motor torque of the i-th propeller, $\boldsymbol{\tau}_{i}$, is,% given by, 
\begin{IEEEeqnarray*}{C}
\boldsymbol{\tau}_{i}=(-1)^{i}b_{T}f_{i}\mathbf{e}_{3}
\end{IEEEeqnarray*}
where the $(-1)^{i}$ term connects each propeller with the correct rotation direction (clockwise and counterclockwise) \cite{Tayebi}.
The total thrust $f$ and moment vector, $^{b}\mathbf{u}$, produced by the propellers are given by,
\begin{IEEEeqnarray}{C}
\begin{bmatrix}
f\\
^b\mathbf{u}
\end{bmatrix}
=
\begin{bmatrix}
1&1&1&1\\
0&d&0&-d\\
-d&0&d&0\\
-b_{T}&b_{T}&-b_{T}&b_{T}\\
\end{bmatrix}\!\!\mathbf{F},
\mathbf{F}=\begin{bmatrix}
f_{1}\\
f_{2}\\
f_{3}\\
f_{4}\\
\end{bmatrix}\IEEEeqnarraynumspace
\label{eq:distribution}
\end{IEEEeqnarray}
with $\mathbf{F}\in\mathbb{R}^{4}$ the thrust vector, and the $4\times4$ matrix to be always full rank for $d,b_{T}\in\mathbb{R}^{+}$.
Thus the thrust magnitude $f$ and the moment $^{b}\mathbf{u}$ will be considered as control inputs with the thrust for each propeller to be calculated from (\ref{eq:distribution}).% for an actual implementation.

The spatial configuration of the quadrotor UAV is described by the quadrotor attitude with respect to the inertial frame and the location of its center of mass, again with respect to the same inertial frame.
The configuration manifold is the special Euclidean group SE(3)=$\mathbb{R}^{3}\times\text{SO(3)}$.
The total thrust produced by the propellers, in the inertial frame, is given by $\mathbf{R}f\mathbf{e}_{3}$.
The equations of motion of the quadrotor system are given by,
\begin{IEEEeqnarray}{rCl}
\mathbf{\dot{x}}&=&{}\mathbf{v}\IEEEnonumber\IEEEeqnarraynumspace\\
m\dot{\mathbf{v}}&=&-mg\mathbf{E}_{3}+\mathbf{R}f\mathbf{e}_{3}\IEEEyesnumber\label{eq:position}\\
\mathbf{J}{}^{b}\dot{\boldsymbol{\omega}}&=&{}^{b}\mathbf{u}-{}^{b}\boldsymbol{\omega}\times\mathbf{J}{}^{b}\boldsymbol{\omega}\IEEEyesnumber\IEEEeqnarraynumspace\label{eq:attitude}\\
\dot{\mathbf{R}}&=&\mathbf{R}S({}^{b}\boldsymbol{\omega}) \IEEEyesnumber
\label{attitude_kinem}\end{IEEEeqnarray}
where $S(.):\mathbb{R}^{3}\rightarrow\mathfrak{so}(3)$ the cross product map given by,
\begin{IEEEeqnarray}{rCl}
\begin{array}{c}
S(\mathbf{r}){=}[{0},{-r_{3}},{r_{2}};{r_{3}},{0},{-r_{1}};{-r_{2}},{r_{1}},0]\\
{S^{-1}}(S(\mathbf{r})){=}\mathbf{r}
\end{array}\label{iso}
\end{IEEEeqnarray}
with the definition of the inverse map, $S^{-1}(.):\mathfrak{so}(3)\rightarrow\mathbb{R}^{3}$, to be given above also.
This identifies the Lie algebra $\mathfrak{so}(3)$ with $\mathbb{R}^{3}$ using the vector cross product in $\mathbb{R}^{3}$.

\section{Quadrotor Tracking Controls}
Given the underactuated nature of quadrotors, in this paper two flight modes are considered:
%Moreover since this paper follows a geometric context/framework it motivates us to introduce two flight modes with each mode tracking a respective/corresponding set of outputs. The flight modes are:
\begin{itemize}
\item Attitude controlled mode:
The controller achieves tracking for the attitude of the quadrotor UAV.
\item Position controlled mode:
The controller achieves tracking for the position vector of the quadrotor CM and a pointing attitude associated with the yaw orientation of the quadrotor UAV.
\end{itemize}
%Using these aforementioned flight modes in suitable successions and by following strictly the switching conditions derived from the rigorous stability proofs the quadrotor caries out a complex desired flight maneuver.
Using these flight modes in suitable successions, a quadrotor can perform a desired complex flight maneuver.
Moreover it will be shown that each mode has stability properties that allow the safe switching between flight modes.
%For example, in this work one of the desired complex flight maneuvers to be tracked consists of a position flight maneuver with the quadrotor translating while gaining altitude, followed by an attitude maneuver were the UAV completes a $360^{o}$ flip, and by a final segment of steady hover at a desired position.
%Any other desired maneuver can be achieved similarly using these two flight modes.

%Our contribution lies in the introduction of new controllers based on nonlinear surfaces with different structure that introduce different dynamics on the closed loop system and with desired control characteristics.
%Furthermore the controllers we introduce allows us to produce a rigorous stability proof and region of attraction without introducing/demanding a restriction on the initial position error like in \cite{geomquadlee} -\cite{qeopidfar}.
%This is an important contribution since that it introduce greater freedom to the user in regards to control objectives.
%Finally we introduce  a lower bound constraint on the controller generated thrusts without interfering with the control objective, ensuring that the generated thrusts stay positive.
%This is a key contribution that differentiate this work from existing geometric controllers \cite{geomquadlee} -\cite{qeo lm,pidfar} since that negative thrusts  can only be produced by a nonstandard quadrotor equipped with variable pitch propellers and aren't realizable in a standard quadrotor configuration.
\newcounter{Prop1}
\addtocounter{Prop1}{1}
\newcounter{Prop2}
\addtocounter{Prop2}{2}
\newcounter{Prop3}
\addtocounter{Prop3}{3}
\newcounter{Prop4}
\addtocounter{Prop4}{4}
\newcounter{Prop5}
\addtocounter{Prop5}{5}
\newcounter{sub}
\addtocounter{sub}{1}
\subsection{Attitude Controlled Mode}
An attitude control system able to follow an arbitrary smooth desired orientation $\mathbf{R}_{d}(t)\in\text{SO(3)}$ and its assosiated angular velocity $^{b}\boldsymbol{\omega}_{d}(t)\in\mathbb{R}^{3}$ is developed next.
\subsubsection{Attitude tracking errors}
For a given tracking command ($\mathbf{R}_{d}$, $^{b}\boldsymbol{\omega}_{d}$) and the current attitude and angular velocity ($\mathbf{R}$, $^{b}\boldsymbol{\omega}$), an attitude error function $\Psi:\text{SO(3)}\times\text{SO(3)}\rightarrow\mathbb{R}$, an attitude error vector $\mathbf{e}_{R}\in\mathbb{R}^{3}$, and an angular velocity error vector $\mathbf{e}_{\omega}\in\mathbb{R}^{3}$ are defined as follows \cite{qeoadapclee}:
\begin{IEEEeqnarray}{rCl}
\Psi(\mathbf{R},\mathbf{R}_{d})&=&\frac{1}{2}tr[\mathbf{I}-\mathbf{R}^{T}_{d}\mathbf{R}]\geq 0
\label{error_function}\\
\mathbf{e}_{R}(\mathbf{R},\mathbf{R}_{d})&=&\frac{1}{2}S^{-1}(\mathbf{R}^{T}_{d}\mathbf{R}-\mathbf{R}^{T}\mathbf{R}_{d})
\label{att_error}\\
\mathbf{e}_{\omega}(\mathbf{R},{}^{b}\boldsymbol{\omega},\mathbf{R}_{d},{}^{b}\boldsymbol{\omega}_{d})&=&{}^{b}\boldsymbol{\omega}-\mathbf{R}^{T}\mathbf{R}_{d}{}^{b}\boldsymbol{\omega}_{d}
\label{ang_vel_error}
\end{IEEEeqnarray}
where $tr[.]$ is the trace function. 
Important properties regarding (\ref{error_function}), (\ref{att_error}), (\ref{ang_vel_error}),
% the attitude error function $\Psi$, the error vectors $\mathbf{e}_{R}, \mathbf{e}_{\omega}$
including the associated attitude error dynamics used throughout this work are included in \text{Proposition \arabic{Prop1}} and \text{Proposition \arabic{Prop2}}  found in Appendix \ref{appA}.
%The desired , is given by (\ref{attitude_kinem}).
\subsubsection{Attitude tracking controller}
A control system is defined for the attitude dynamics (attitude controlled mode) of the quadrotor UAV stabilizing $\mathbf{e}_{R}$, $\mathbf{e}_{\omega}$ to zero exponentially.

\textbf{Proposition \arabic{Prop3}.} 
For $\eta,k_{R},k_{\omega}\in\mathbb{R}^{+}$, with initial conditions satisfying,
%For an arbitrary smooth desired orientation $\mathbf{R}_{d}(t)\in\text{SO(3)}$, an angular velocity $^{b}\boldsymbol{\omega}_{d}(t)\in\mathbb{R}^{3}$ and initial conditions, 
\begin{IEEEeqnarray}{C}
\Psi(\mathbf{R}(0),\mathbf{R}_{d}(0))<2
\label{Psi_0}\\
\lVert\mathbf{e}_{\omega}(0)\rVert^{2}<2\eta k_{R}\left(2-\Psi(\mathbf{R}(0),\mathbf{R}_{d}(0))\right)\label{surface_0}
\end{IEEEeqnarray}
%with (\ref{Psi_0}-\ref{surface_0}) referring to the initial conditions of the selected nonlinear surface,
and for a desired arbitrary smooth attitude $\mathbf{R}_{d}(t)\in\text{SO(3)}$ in,
\begin{IEEEeqnarray}{rCl}
L_{2}&=&\{(\mathbf{R},\mathbf{R}_{d})\in\text{SO(3)}\times\text{SO(3)}|\Psi(\mathbf{R},\mathbf{R}_{d})<2 \}
\label{L_2}
\end{IEEEeqnarray}
%i.e., $\mathbf{R}_{d}(t)$ is not antipodal to $\mathbf{R}(t)$.
then, under the assumption of perfect parameter knowledge, we propose 
%for the attitude tracking of the quadrotor UAV, 
the following nonlinear surface-based controller,
\begin{IEEEeqnarray}{rCl}
^{b}\mathbf{u}&=&{}^{b}\boldsymbol{\omega}\times\mathbf{J}{}^{b}\boldsymbol{\omega}-\mathbf{J}\left(\frac{k_{R}}{k_{\omega}}\dot{\mathbf{e}}_{R}+\mathbf{a}_{d}+\eta\mathbf{s}_{R}\right)\label{att_contr}
\end{IEEEeqnarray}
where $\mathbf{a}_{d}$ is defined in (\ref{E_ad}) and the nonlinear surface $\mathbf{s}_{R}$ is defined by,
\begin{IEEEeqnarray}{C}
\mathbf{s}_{R}=k_{R}\mathbf{e}_{R}+k_{\omega}\mathbf{e}_{\omega}\label{att_surface}
\end{IEEEeqnarray}
then 
%under the assumption of perfect parameter knowledge,
the zero equilibrium of the quadrotor closed loop attitude tracking error $(\mathbf{e}_{R},\mathbf{e}_{\omega})=(\mathbf{0},\mathbf{0})$ is almost globally exponentially stable.
Moreover there exist constants $\mu,\tau>0$ such that
\begin{IEEEeqnarray}{C}
\Psi(\mathbf{R},\mathbf{R}_{d})<min\{2,\mu e^{-\tau t}\}
\label{Psi_bou}
\end{IEEEeqnarray}
%\textbf{Proof.}

\textbf{Proof of Proposition \arabic{Prop3}.}
We employ a sliding methodology in $L_{2}$ by defining the nonlinear surface in terms of the attitude configuration errors (\ref{att_error}), (\ref{ang_vel_error}) and apply Lyapunov analysis.
\begin{enumerate}[(a)]
\item Lyapunov candidate: We define,
\begin{IEEEeqnarray}{rCl}
V&=&\frac{1}{2k_{\omega}}\mathbf{s}_{R}^{T}\mathbf{s}_{R}+2\eta k_{R}k_{\omega}\Psi
\label{att_lyap}
\end{IEEEeqnarray}
Differentiating (\ref{att_lyap}) and substituting (\ref{att_contr}) we get,
\begin{IEEEeqnarray}{rCl}
\dot{V}&=&-\eta\mathbf{z}^{T}_{R}\mathbf{W}_{3}\mathbf{z}_{R},\mathbf{W}_{3}=\begin{bmatrix}
k_{R}^{2}&0\\
0&k_{\omega}^{2}
\end{bmatrix}
\label{Datt_lyap}
\end{IEEEeqnarray}
where $\mathbf{z}_{R}=[\lVert\mathbf{e}_{R}\rVert;\lVert\mathbf{e}_{\omega}\rVert]$.
\item Boundedness of $\Psi$: We define the Lyapunov function,
\begin{IEEEeqnarray}{rCl}
V_{\Psi}&=&\frac{1}{2}\mathbf{e}_{\omega}^{T}\mathbf{e}_{\omega}+\eta k_{R}\Psi\label{ater_lyap}\\
%\dot{V}_{\Psi}&\leq&-(\eta k_{\omega}+\frac{k_{R}}{k_{\omega}})\lVert\mathbf{e}_{\omega}\rVert^{2}\leq0
\dot{V}_{\Psi}&\leq&-\eta k_{\omega}\lVert\mathbf{e}_{\omega}\rVert^{2}\leq0
\label{dater_lyap}
\end{IEEEeqnarray}
 Equations (\ref{ater_lyap}-\ref{dater_lyap}) imply that $V_{\Psi}(t)\leq V_{\Psi}(0),\forall t\geq 0$.
Applying (\ref{surface_0}) we obtain,
\begin{IEEEeqnarray}{C}
\eta k_{R} \Psi(\mathbf{R}(t),\mathbf{R}_{d}(t)){\leq} V_{\Psi}(t){\leq} V_{\Psi}(0){<}2\eta k_{R} \label{hkrko_Psi}
\end{IEEEeqnarray}
Implying that the attitude error function is bounded by,
\begin{IEEEeqnarray}{C}
\Psi(\mathbf{R}(t),\mathbf{R}_{d}(t))\leq \psi_{a} < 2,\forall t\geq 0\label{B_Psi}
\end{IEEEeqnarray}
where $\psi_{a}={V(0)}/{\eta k_{R} }$. Thus $\mathbf{R}(t)\in L_{2}$.
\item Exponential Stability: Using (\ref{low_Psi}), (\ref{upper_Psi}) it follows that $V$ is bounded,
\begin{IEEEeqnarray}{C}
\mathbf{z}^{T}_{R}\mathbf{W}_{1}\mathbf{z}_{R}\leq V\leq\mathbf{z}^{T}_{R}\mathbf{W}_{2}\mathbf{z}_{R}
\end{IEEEeqnarray}
where $\mathbf{W}_{1}$, $\mathbf{W}_{2}$ are positive definite matrices given by,
\begin{IEEEeqnarray}{C}
\mathbf{W}_{1}=
\begin{bmatrix}
\frac{k_{R}^{2}}{2k_{\omega}}+\eta k_{R}k_{\omega}&-\frac{k_{R}}{2}\\
-\frac{k_{R}}{2}&\frac{k_{\omega}}{2}
\end{bmatrix},\\
\mathbf{W}_{2}=
\begin{bmatrix}
\frac{k_{R}^{2}}{2k_{\omega}}+\frac{2}{2-\psi_{a}}\eta k_{R}k_{\omega}&\frac{k_{R}}{2}\\
\frac{k_{R}}{2}&\frac{k_{\omega}}{2}
\end{bmatrix}
\end{IEEEeqnarray}
Thus the following inequalities hold,
\begin{IEEEeqnarray}{C}
\lambda_{min}(\mathbf{W}_{1})\lVert \mathbf{z}_{R} \rVert^{2}\leq V \leq\lambda_{max}(\mathbf{W}_{2})\lVert\mathbf{z}_{R}\rVert^{2}\\
\dot{V} \leq -\eta\lambda_{min}(\mathbf{W}_{3})\lVert\mathbf{z}_{R}\rVert^{2}
\end{IEEEeqnarray}
Then for $\tau=\frac{\eta\lambda_{min}(\mathbf{W}_{3})}{\lambda_{max}(\mathbf{W}_{2})}$ the following holds,
\begin{IEEEeqnarray}{C}
\dot{V} \leq -\tau V
\end{IEEEeqnarray}
Thus the zero equilibrium of the attitude tracking error $\mathbf{e}_{R}$, $\mathbf{e}_{\omega}$ is exponentially stable.
Using (\ref{upper_Psi}) then,
\begin{IEEEeqnarray}{C}
%(2-\psi_{a})\lambda_{min}(\mathbf{W}_{1})\Psi \leq \lambda_{min}(\mathbf{W}_{1})\lVert \mathbf{e}_{R} \rVert^{2}\leq \lambda_{min}(\mathbf{W}_{1})\lVert \mathbf{z}_{R} \rVert^{2}\leq V(t) \leq V(0)e^{-\tau t}\IEEEnonumber
(2-\psi_{a})\lambda_{min}(\mathbf{W}_{1})\Psi \leq  V(t) \leq V(0)e^{-\tau t}\IEEEnonumber
\end{IEEEeqnarray}
Thus $\Psi$ exponentially decreases and from (\ref{B_Psi}) we arrive to (\ref{Psi_bou}).
%This completes the proof.
This completes the proof. $\blacksquare$
\end{enumerate}
%See Appendix \ref{appA}

The region of attraction given by (\ref{Psi_0})-(\ref{surface_0}) ensures that the initial attitude error is less than $180^{o}$  with respect to an axis-angle rotation for a desired $\mathbf{R}_{d}$ (i.e., $\mathbf{R}_{d}(t)$ is not antipodal to $\mathbf{R}(t)$).
This is the best that one can do since (\ref{att_error}) vanishes at the antipodal equilibrium.
Consequently exponential stability is guaranteed almost globally (everywhere except the antipodal equilibrium).
This was expected since in \cite{obstruction}, it was shown that the topology of SO(3) prohibits the design of a smooth global controller.
Since (\ref{att_contr}) is developed directly on SO(3), it completely avoids singularities and ambiguities associated with minimum attitude representations like Euler angles or quaternions.
Also the controller can be applied to the attitude dynamics of any rigid body and not only on quadrotor systems.

Note that, tracking of the quadrotor attitude does not require the definition of the thrust magnitude $f$.
%f=\sum_{i=1}^{4} f_i$.
Therefore no explicit regulation of the quadrotor position takes place.
Thus the attitude mode is better suited for short durations of time.
In \cite{geomquadlee},\cite{geommac},\cite{geomquadlee_asian}, (\ref{eq:distribution}) is used in conjunction with a suitable expression for $f$, to track a desired altitude command.
A different approach in regards to this matter is taken here and is presented in Section \ref{overA}.
\subsection{Position Controlled Mode}
%A control system is defined for the position dynamics (position controlled mode) of the quadrotor UAV stabilizing $\mathbf{e}_{x}$, $\mathbf{e}_{v}$, $\mathbf{e}_{R}$, $\mathbf{e}_{\omega}$ to zero asymptotically while for initial conditions in the sublevel set (\ref{L_2}) the proposed controller renders the zero equilibrium exponentially stable.
\subsubsection{Position tracking errors}
For an arbitrary smooth position tracking instruction $\mathbf{x}_{d}\in\mathbb{R}^{3}$, the tracking errors for the position and the velocity are taken as,
\begin{IEEEeqnarray}{C}
\mathbf{e}_{x}=\mathbf{x}-\mathbf{x}_{d},\; \mathbf{e}_{v}=\mathbf{v}-\dot{\mathbf{x}}_{d}\label{pos_error}
\end{IEEEeqnarray}
For $a,k_{x},k_{v}{\in}\mathbb{R}^{+}$ the position nonlinear surface is defined as,
\begin{IEEEeqnarray}{C}
\mathbf{s}_{x}=k_{x}\mathbf{e}_{x}+k_{v}\mathbf{e}_{v}\label{pos_surf}
\end{IEEEeqnarray}

%In this mode the attitude dynamics must be synergetic with the desired position tracking instruction.
In this mode the attitude dynamics must be compatible with the desired position tracking instruction.
This results in the definition of a position-induced attitude and using it as a command for the attitude dynamics.
Thus the attitude dynamics are guided to follow the position related attitude $\mathbf{R}_{x}(t)\in\text{SO(3)}$ and angular velocity $^{b}\boldsymbol{\omega}_{x}(t)$ given by,
\begin{IEEEeqnarray}{C}
\mathbf{R}_{x}{=}\left[\mathbf{e}_{1_{h}},\frac{\mathbf{e}_{3_{x}}\times\mathbf{e}_{1_{h}}}{\lVert\mathbf{e}_{3_{x}}\times\mathbf{e}_{1_{h}}\rVert},\mathbf{e}_{3_{x}}\right],\;{}^{b}\boldsymbol{\omega}_{x}{=}{S^{-1}}\!(\mathbf{R}^{T}_{x}\dot{\mathbf{R}}_{x})\label{Rxbox}\\
\mathbf{e}_{3_{x}}{=}\frac{mg\mathbf{E}_{3}-m\frac{k_{x}}{k_{v}}\mathbf{e}_{v}-a\mathbf{s}_{x}+m\ddot{\mathbf{x}}_{d}}{\lVert{mg\mathbf{E}_{3}-m\frac{k_{x}}{k_{v}}\mathbf{e}_{v}-a\mathbf{s}_{x}+m\ddot{\mathbf{x}}_{d}}\rVert}\in\text{S}^{2}\label{e_3_x}
\end{IEEEeqnarray}
and by trajectory design the denominator of (\ref{e_3_x}) is non zero.

The vector $\mathbf{e}_{1_{h}}$ is chosen as in \cite{geommac}, by defining arbitrarily a desired direction $\mathbf{e}_{1_{d}}\in\text{S}^{2}$ of the first body fixed axis of the quadrotor such that $\mathbf{e}_{1_{d}}\nparallel\mathbf{e}_{3_{x}}$.
The heading is then given by 
%$\mathbf{e}_{1_{h}}=\frac{(\mathbf{e}_{3_{x}}\times\mathbf{e}_{1_{d}})\times\mathbf{e}_{3_{x}}}{\lVert(\mathbf{e}_{3_{x}}\times\mathbf{e}_{1_{d}})\times\mathbf{e}_{3_{x}}\rVert}$
\begin{IEEEeqnarray*}{C}
\mathbf{e}_{1_{h}}=\frac{(\mathbf{e}_{3_{x}}\times\mathbf{e}_{1_{d}})\times\mathbf{e}_{3_{x}}}{\lVert(\mathbf{e}_{3_{x}}\times\mathbf{e}_{1_{d}})\times\mathbf{e}_{3_{x}}\rVert}
\end{IEEEeqnarray*}
and $\mathbf{e}_{1_{h}}\perp\mathbf{e}_{3_{x}}$, thus $\mathbf{R}_{x}(t)\in\text{SO(3)}$. 
%Thus the generated $\mathbf{R}_{x}\in\text{SO(3)}$.
\subsubsection{Position tracking controller}
A control system is developed for the position dynamics (position controlled mode) of the quadrotor UAV achieving almost global asymptotic stabilization of $(\mathbf{e}_{x},\mathbf{e}_{v},\mathbf{e}_{R},\mathbf{e}_{\omega})$ to the zero equilibrium.

For a sufficiently smooth pointing direction $\mathbf{e}_{1_{d}}(t)\in\text{S}^{2}$, and position tracking instruction $\mathbf{x}_{d}(t)\in\mathbb{R}^{3}$ the following position controller is defined, 
%composed by (\ref{att_contr}) and by the thrust magnitude,
\begin{IEEEeqnarray}{rCl}
f(\mathbf{x}_{d},\dot{\mathbf{x}}_{d},\ddot{\mathbf{x}}_{d})&{=}&(mg\mathbf{E}_{3}{-}m\frac{k_{x}}{k_{v}}\mathbf{e}_{v}{-}a\mathbf{s}_{x}{+}m\ddot{\mathbf{x}}_{d})^{T}\mathbf{R}\mathbf{e}_{3}\IEEEyessubnumber\label{f}\\
^{b}\mathbf{u}(\mathbf{R}_{x},{{}^{b}\boldsymbol{\omega}_{x}})&{=}&{}^{b}\boldsymbol{\omega}{\times}\mathbf{J}{}^{b}\boldsymbol{\omega}{-}\mathbf{J}\left(\frac{k_{R}}{k_{\omega}}\dot{\mathbf{e}}_{R_{x}}{+}\mathbf{a}_{d_{x}}{+}\eta\mathbf{s}_{R_{x}}\right)\IEEEyessubnumber\label{att_contrx}
\end{IEEEeqnarray}
where $\dot{\mathbf{e}}_{R_{x}},\mathbf{a}_{d_{x}},\mathbf{s}_{R_{x}}$ are given by (\ref{dot_Att_Error}), (\ref{E_ad}), (\ref{att_surface}), through the use of (\ref{Rxbox}). 
The closed loop system is  defined by (\ref{eq:position})-(\ref{attitude_kinem}) under the action of (\ref{f})-(\ref{att_contrx}).
We proceed to state the result of exponential stability
%that zero
%$(\mathbf{e}_{x},\mathbf{e}_{v},\mathbf{e}_{R},\mathbf{e}_{\omega})=(\mathbf{0},\mathbf{0},\mathbf{0},\mathbf{0})$ 
%is an exponentially stable equilibrium 
of the quadrotor closed loop dynamics.

\textbf{Proposition \arabic{Prop4}.}
%For a sufficiently smooth pointing direction $\mathbf{e}_{1_{d}}(t)\in\text{S}^{2}$, and position tracking instruction $\mathbf{x}_{d}(t)\in\mathbb{R}^{3}$ the following position controller is defined, 
%composed by (\ref{att_contr}) and by the thrust magnitude,
%\begin{IEEEeqnarray}{rCl}
%f&=&(mg\mathbf{E}_{3}-m\frac{k_{x}}{k_{v}}\mathbf{e}_{v}-a\mathbf{s}_{x}+m\ddot{\mathbf{x}}_{d})\cdot\mathbf{R}\mathbf{e}_{3}\label{f}\\
%^{b}\mathbf{u}&=&{}^{b}\boldsymbol{\omega}\times\mathbf{J}{}^{b}\boldsymbol{\omega}-\mathbf{J}\left(\frac{k_{R}}{k_{\omega}}\dot{\mathbf{e}}_{R}+\mathbf{a}_{d}+\eta\mathbf{s}\right)\label{att_contr}
%\end{IEEEeqnarray}
For initial conditions in the domain,
%For $\eta,k_{R},k_{\omega}\in\mathbb{R}^{+}$, initial conditions,
\begin{IEEEeqnarray}{rCl}
D_{x}&=&\{(\mathbf{e}_{x},\mathbf{e}_{v},\mathbf{e}_{R},\mathbf{e}_{\omega})\in\mathbb{R}^{3}\times\mathbb{R}^{3}\times\mathbb{R}^{3}\times\mathbb{R}^{3}|\IEEEnonumber\\
&{ }&\Psi(\mathbf{R}(0),\mathbf{R}_{x}(0))<\psi_{p}<1\}
\label{D_x}
\end{IEEEeqnarray}
%\begin{IEEEeqnarray}{C}
%\Psi(\mathbf{R}(0),\mathbf{R}_{x}(0))<\psi_{p}<1,\lVert\mathbf{e}_{r}\rVert<\theta\label{pos_psi}
%\end{IEEEeqnarray}
a uniformly bounded desired acceleration, 
\begin{IEEEeqnarray}{C}
\lVert mg\mathbf{E}_{3}+m\ddot{\mathbf{x}}_{d}\rVert \leq B \label{B}
\end{IEEEeqnarray}
with $B\in\mathbb{R}^{+}$.
We define $\mathbf{\Pi}_{1},\mathbf{\Pi}_{2}\in\mathbb{R}^{2\times2}$ as,
\begin{IEEEeqnarray}{C}
\mathbf{\Pi}_{1}{=}
\begin{bmatrix}
ak_{x}^{2}(1{-}\theta)&-ak_{x}k_{v}\theta{-}\frac{mk_{x}^{2}\theta}{2k_{v}}\\
-ak_{x}k_{v}\theta{-}\frac{mk_{x}^{2}\theta}{2k_{v}}&ak_{v}^{2}{-}\theta(mk_{x}{+}ak_{v}^{2})
\end{bmatrix},\IEEEnonumber\\
\mathbf{\Pi}_{2}=
\begin{bmatrix}
Bk_{x}&0\\
Bk_{v}&0
\end{bmatrix}\label{P}
\end{IEEEeqnarray}
where $\theta<\theta_{max}\in\mathbb{R}^{+}$.
The attitude error bound, $\psi_{p}$, satisfies, 
%$\theta_{max}=\sqrt{\psi_{p}(2-\psi_{p})}$
\begin{IEEEeqnarray*}{C}
\theta_{max}=\sqrt{\psi_{p}(2-\psi_{p})}
\end{IEEEeqnarray*}
and $\theta_{max}$ is given by,
\begin{IEEEeqnarray}{rCl}
\theta_{max}&=&\min\{\frac{ak_{v}^{2}}{ak_{v}^{2}{+}mk_{x}},\delta_{1}+\delta_{2}\},\IEEEyesnumber\label{theta}\\
\delta_{1}&=&2{\frac {k_{v}^{2}\sqrt {4k_{x}^{4}k_{v}^{4}a^{4}+4k_{x}^{5}k_{v}^{2}{a
}^{3}m+2k_{x}^{6}{m}^{2}{a}^{2}}}{k_{x}^{4}{m}^{2}}}
\IEEEnonumber\\
\delta_{2}&=&-4{\frac {{a}^{2}k_{v}^{4}}{{m}^{2}k_{x}^{2}}}{-}2{\frac {ak_{v}^{2}}{mk_{x}}}\IEEEnonumber
\end{IEEEeqnarray}
%The attitude error bound, $\psi_{p}$, satisfies $\theta_{max}=\sqrt{\psi_{p}(2-\psi_{p})}$.
In conjunction with suitable gains $\eta,k_{R},k_{\omega}\in\mathbb{R}^{+}$, such that,
\begin{IEEEeqnarray}{C}
\lambda_{min}(\mathbf{W}_{3})>\frac{\lVert\mathbf{\Pi}_{2}\rVert^{2}}{4\eta\lambda_{min}(\mathbf{\Pi}_{1})}\label{w3}
\end{IEEEeqnarray}
%where $\mathbf{W}_{3}$ is defined in (\ref{Datt_lyap}) and 
then the zero equilibrium of the closed loop errors $(\mathbf{e}_{x},\mathbf{e}_{v},\mathbf{e}_{R},\mathbf{e}_{\omega})$ is exponentially stable. A region of attraction is identified by (\ref{D_x}), and
\begin{IEEEeqnarray}{C}
\lVert\mathbf{e}_{\omega}(0)\rVert^{2}<2\eta k_{R}\left(\psi_{p}-\Psi(\mathbf{R}(0),\mathbf{R}_{x}(0))\right)\label{ep_0}
\end{IEEEeqnarray}
\textbf{Proof.}
See Appendix \ref{appB}.

Proposition \arabic{Prop4} requires that the norm of the initial attitude error is less than % $\lVert\mathbf{e}_{R}\rVert<\theta_{max}$
$\theta_{max}$ to achieve exponential stability (the upper bound of theta, (\ref{theta}), depends solely on the control gains and the quadrotor mass) which corresponds to a slightly reduced region of attraction in comparison to the regions in \cite{geomquadlee} -\cite{qeopidfar}.
This is because no restriction on the initial position/velocity error was applied during the stability proof.
This approach is not only novel (wrt. the geometric literature) but it also offers the advantage of simplifying the trajectory design procedure, in regards to designer considerations.
% (simple criteria are provided).
In contrast to the proposed approach the region of attraction proposed in other geometric treatments includes bounds on the initial position or velocity (see \cite{geomquadlee} -\cite{qeopidfar}) meaning that the trajectory designer should comply to the position/velocity bounds and also to the attitude bound, a more involved/complicated task.
%To our knowledge a controller with a position/velocity free region of attraction is a new contribution.

If a user prefers a larger basin of attraction, this can be achieved by introducing bounds on the initial position/velocity.
Then two new regions of attraction are produced involving larger initial attitude errors and are given by (\ref{B}), (\ref{ep_0}) and,
\begin{IEEEeqnarray}{C}
\Psi(\mathbf{R}(0),\mathbf{R}_{x}(0))<\psi_{p}<1, \lVert\mathbf{e}_{x/v}(0)\rVert<e_{x/v_{max}}\label{fra_the_xv}\\
\theta<\theta_{max}=\min\{\frac{ak_{v}^{2}}{ak_{v}^{2}{+}mk_{x}}\}
%,\theta_{max}=\sqrt{\psi_{p}(2-\psi_{p})}
\IEEEyesnumber\label{thetaxv}
\end{IEEEeqnarray}%$\lVert\mathbf{e}_{x/v}(0)\rVert{<}e_{x/v_{max}}$
where the second inequality in (\ref{fra_the_xv}) denotes either a bound on the initial position error, $e_{x_{max}}$, or a bound on the initial velocity error, $e_{v_{max}}$, but not on both (see Appendix \ref{appB} for more details and expressions regarding $\Pi_{1}$, $\Pi_{2}$, that comply with (\ref{w3})).
Depending on the user preference the trajectory design procedure can be realized using either one of the three regions of attraction guiding us to favorable conditions for switching between flight modes.
For completeness all three regions of exponential stability were derived, but this work focuses on the position/velocity free basin.

Finally, the proposition that follows shows that the structure of the position controller is characterized by almost global exponential attractiveness.
This compensates for the reduced position/velocity free region of attraction and introduces greater freedom to the user in regards to control objectives since the region of attraction does not depend explicitly on the initial position/velocity error.
If the quadrotor initial states are outside of $D_{x}$ in (\ref{D_x}), Proposition \arabic{Prop3} still applies.
Thus the attitude state enters in (\ref{D_x}) at a finite time $t^{*}$ and the results of Proposition \arabic{Prop4} take effect.
The result regarding the position controlled mode is stated next.

\textbf{Proposition \arabic{Prop5}.} For initial conditions satisfying (\ref{surface_0}), and
\begin{IEEEeqnarray}{C}
\psi_{p}\leq\Psi(\mathbf{R}(0),\mathbf{R}_{x}(0))<2\label{Psi_atr}
\end{IEEEeqnarray}
and a uniformly bounded desired acceleration (\ref{B}), the thrust magnitude defined in (\ref{f}), in conjunction with the control moment (\ref{att_contrx}), renders the zero equilibrium of $(\mathbf{e}_{x},\mathbf{e}_{v},\mathbf{e}_{R},\mathbf{e}_{\omega})$ almost globally exponentially attractive.

\textbf{Proof of Proposition \arabic{Prop5}.}
See Proposition 4 in \cite{geommac} but apply the proposed thrust (\ref{f}), to arrive to the same conclusion.

%The procedure is summarized next. 
%A positive definite function is defined as a function of the position errors,
%\begin{IEEEeqnarray}{C}
%V_{a}=\frac{\mathbf{e}_{x}^{T}\mathbf{e}_{x}}{2}+\frac{m\mathbf{e}_{v}^{T}\mathbf{e}_{v}}{2}\label{Va}
%\end{IEEEeqnarray}
%Then it holds that $\lVert\mathbf{e}_{x}\rVert\leq\sqrt{2V_{a}}$, $\lVert\mathbf{e}_{v}\rVert\leq\sqrt{({2}/{m})V_{a}}$. We differentiate (\ref{Va}) and use equations (\ref{eq:position}), (\ref{pos_error}), (\ref{B}), (\ref{f}) to get,
%\begin{IEEEeqnarray}{C}
%\dot{V}_{a}\leq d_{1}V_{a}+d_{2}\sqrt{V_{a}}\IEEEnonumber\\
%d_{1}=2\{(1+ak_{x})\sqrt{\frac{1}{m}}+(\frac{k_{x}}{k_{v}}+\frac{ak_{v}}{m})\},d_{2}=2B\sqrt{\frac{2}{m}}\IEEEnonumber
%\end{IEEEeqnarray}
%Suppose $V_{a}\geq 1$ for a time interval $[t_{a},t_{b}]\subset[0,t^{*}]$, then $\sqrt{V_{a}}\leq V_{a}$ and the following holds,
%\begin{IEEEeqnarray}{C}
%\dot{V}_{a}\leq (d_{1}+d_{2})V_{a}\rightarrow V_{a}(t)\leq V_{a}(t_{a})e^{(d_{1}+d_{2})(t-t_{a})}\IEEEnonumber
%\end{IEEEeqnarray}
%This implies that $V_{a}$ is bounded for $0\leq t\leq t^{*}$.

In other words, Proposition \arabic{Prop5} shows that during the finite time that it takes for the attitude states to enter the region of attraction for exponential stability (\ref{D_x}), (\ref{ep_0}), the position tracking error $\lVert\mathbf{z}_{x}\rVert$ remains bounded.
The calculated region of exponential attractiveness given by (\ref{Psi_atr}) ensures that the initial attitude error is less than $180^{o}$  with respect to an axis-angle rotation for a desired $\mathbf{R}_{x}$ (i.e., $\mathbf{R}_{x}(t)$ is not antipodal to $\mathbf{R}(t)$).
Consequently the zero equilibrium of the tracking errors is almost globally exponentially attractive.
%\section{Adaptive Quadrotor Tracking Controls}
%\subsection{Attitude Controlled Mode}
%\subsection{Position Controlled Mode}
\section{Exploiting Overactuation\label{overA}}
As mentioned earlier, during aggressive maneuvers, existing geometric controllers produce negative thrusts \cite{geomquadlee},\cite{geommac},\cite{geomquadlee_asian},\cite{qeopidfar}, that are not realizable with standard quadrotors.
When the desired thrust is negative, in a standard quadrotor, the controller drives the propeller speed to zero (a saturation state) in an attempt to achieve the thrust.
This action has two adverse effects.
Firstly and most obviously the tracking error increases significantly since the desired control effort is not available for the maneuver and secondly the out-runner motors undergo an aggressive state change where they need to come to a complete halt and again instantaneously achieve a high RPM count.
This is not only strenuous for the motors and reduces their lifespan, it also is extremely expensive energy-wise reducing the available flight time of the UAV.
Another important consideration is that the stability proofs accompanying the controllers do not account for thrust saturations and this also holds for \cite{geomquadlee},\cite{geommac},\cite{geomquadlee_asian},\cite{qeopidfar} and the majority of geometric controllers in the bibliography.
Thus in order for the stability proofs (regions of attraction) to hold the desired control effort must be available, i.e., avoid saturation or/and change of sign.
%Also it is very important to emphasize that the guarantees produced by the rigorous stability proofs hold since the controller generated thrusts produced by the GNCS are generated in a non invasive manner through (\ref{xi}) and comply to the controller demands.
%This is the most important point of this analysis, since the stability proof does not account for thrust saturations and this also holds for \cite{geomquadlee},\cite{geommac},\cite{qeopidfar}.
%Thus it is very important that the desired control effort is available, i.e., avoid saturation.
%However, negative thrusts  can only be generated by a nonstandard quadrotor equipped with variable pitch propellers (or other means) and are not realizable in a standard quadrotor configuration.

By studying the occurrence of negative thrusts through extensive simulations, it was observed that thrusts remain positive if the control task at hand is a position trajectory of a relatively reasonable rate.
On the other hand, if the control task entails a large angle attitude maneuver the thrusts can certainly become negative, even if the attitude maneuver is conducted at very slow rates.
Therefore, it is important to develop a method, realizable in real time, to distribute the generated control effort (\ref{att_contr}), to the four thrusters of the quadrotor without interfering with the control objective while simultaneously complying with the following constraint,
%The above observation motivated us to research this important phenomenon and try to find a realizable/ effective solution.
% to this constrained optimization problem.
% The exact problem formulation is given in the following proposition.
%\textbf{Problem.} Find a method, that is realizable in real-time, to distribute the generated control effort (\ref{att_contr}), to the four thrusters of the quadrotor without interfering with the control objective while simultaneously complying with the following constraint,
\begin{IEEEeqnarray}{C}
f_{max}>f_{i}>f_{min},i=1,..,4\label{f_con}%\\
%f_{b}>0,f_{b}\in\mathbb{R}\IEEEnonumber
\end{IEEEeqnarray}
%where $f_{b}\in\mathbb{R}$ is a positive constant.

%The problem above poses as a constrained optimization problem with a plethora of solutions in the bibliography \cite{Boyd}.
%This poses as a constrained optimization problem for which a number of solutions exists \cite{Boyd}. %[angeles,zexiang]
This poses as a constrained optimization problem for which many solutions exists \cite{Boyd}. %[angeles,zexiang]
However, to take advantage of the dynamics of the system, in conjunction with the proposed control strategy, an alternative solution is proposed that is extremely simple, fast and complies with the conditions above.

The solution starts with the realization that even though the quadrotor is underactuated in SE(3), it can be viewed as an overactuated platform in SO(3), the configuration space of its attitude dynamics. %, since that it is equipped with four thrusters.
Resultantly it was identified that during the Attitude Controlled Mode, this aforementioned actuation redundancy allows to achieve (\ref{att_contr}) and additional constraints.
The moment vector, $^{b}\mathbf{u}$, to be produced by the propellers is associated with the thrust vector, $\mathbf{F}\in\mathbb{R}^{4}$, by,
\begin{IEEEeqnarray}{rCl}
\mathbf{F}=
%\begin{bmatrix}
%f_{1}\\
%f_{2}\\
%f_{3}\\
%f_{4}\\
%\end{bmatrix}=
\mathbf{A}^{\#}({}^b\mathbf{u})&,&\mathbf{A}=
\begin{bmatrix}
0&d&0&-d\\
-d&0&d&0\\
-b_{T}&b_{T}&-b_{T}&b_{T}\\
\end{bmatrix}\IEEEeqnarraynumspace\label{alloc}\\
\mathbf{A}^{\#}&=&\mathbf{A}^{T}(\mathbf{A}\mathbf{A}^{T})^{-1}\in\mathbb{R}^{4\times3}
\label{pinv}
\end{IEEEeqnarray}
with $\mathbf{A}$ to always have full row rank when $d,b_{T}\neq0$ and $\mathbf{A}^{\#}$ is the Moore-Penrose pseudoinverse.
%The null space of (\ref{alloc}) is exploited (see \cite{Siciliano}, ch. 3.5.1) to achieve additional tasks,
The null space of (\ref{alloc}) is exploited to achieve additional tasks, \cite{Siciliano},
\begin{IEEEeqnarray}{L}
%[f_{1},f_{2},f_{3},f_{4}]^{T}
\mathbf{F}{=}\mathbf{A}^{\#}{}^b\mathbf{u}{+}(\mathbf{I}{-}\mathbf{A}^{\#}\mathbf{A})\boldsymbol{\xi}\IEEEyesnumber
%\boldsymbol{\xi}{=}{-}k_{H}\nabla H
\label{nullspa}
\end{IEEEeqnarray}
where $\boldsymbol{\xi}\in\mathbb{R}^{4}$ is a suitable vector designed to achieve two objectives.
Namely avoiding saturations as a first priority and secondly allow for the quadrotor during the attitude maneuver to track a desired position.
%where for $\iota_{x},\iota_{y},\iota_{z}\in\mathbb{R}^{+},\boldsymbol{\xi}$ is a suitable vector given by,
The designed vector $\boldsymbol{\xi}$ is given by,
\begin{IEEEeqnarray}{L}
\boldsymbol{\xi}{=}\int^{t2}_{t1}\!\!{\nabla}_{\mathbf{F}}H(\mathbf{F})d\tau{+}{\begin{bmatrix}
1&1&1&1\\
0&d&0&-d\\
-d&0&d&0\\
-b_{T}&b_{T}&-b_{T}&b_{T}\\
\end{bmatrix}}^{\!-1}
\!{\begin{bmatrix}
f_{p}\\%f_{alt}
0\\
0\\
0\\
\end{bmatrix}}\IEEEyesnumber
\label{xi}\\
%f_{alt}=U_{3}(\mathbf{e}_{3}^{T}(\mathbf{R}\mathbf{e}_{3}))\IEEEyesnumber
%f_{p}{=}\Big({\begin{bmatrix}
%\iota_{x}&0&0\\
%0&\iota_{y}&0\\
%0&0&\iota_{z}\\
%\end{bmatrix}}{(mg\mathbf{E}_{3}{-}m\frac{k_{x}}{k_{v}}\mathbf{e}_{v})}\Big){\cdot}\mathbf{R}\mathbf{e}_{3}
f_{p}{=}{\left({\begin{bmatrix}
\iota_{x}&0&0\\
0&\iota_{y}&0\\
0&0&\iota_{z}\\
\end{bmatrix}}{(mg\mathbf{E}_{3}{-}m\frac{k_{x}}{k_{v}}\mathbf{e}_{v}{-}k_{\xi}\mathbf{s}_{x}{+}m\ddot{\mathbf{x}}_{d})}\right)}^{\!\!\!T}\mathbf{R}\mathbf{e}_{3}\IEEEnonumber
\end{IEEEeqnarray}

The first component of (\ref{xi}) is responsible for applying actuator constraints by keeping $f_{i}$ as close to $f_{idl}$ and between $f_{min}$ and $f_{max}$ through the gradient of a suitable function, $H(\mathbf{F})$.
The function $H(\mathbf{F})=\sum_{i=1}^{4}h(f_{i})$ is chosen as in \cite{Ryll}, 
%where $U_{3}$ is the third component of the vector given in (\ref{U}) and the function $H(\mathbf{F})=\sum_{i=1}^{4}h(f_{i})$ is chosen as in \cite{Ryll},
%\begin{IEEEeqnarray*}{L}
%h(f_{i})=
%\begin{cases}
%k_{h_{1}}tan^{2}(\gamma_{1}\lvert f_{i}\rvert+\gamma_{2}), & f_{min}{<}\lvert f_{i}\rvert{\leq}f_{idl}\\
%\frac{k_{h_{2}}}{2}(\lvert f_{i}\rvert-f_{idl})^{2}+\frac{(\lvert f_{i}\rvert-f_{idl})^{2}}{(\lvert f_{i}\rvert-f_{max})},& \lvert f_{i}\rvert>f_{idl}
%\end{cases}
%\end{IEEEeqnarray*}
\begin{IEEEeqnarray*}{L}
h(f_{i})=
\begin{cases}
k_{h_{1}}tan^{2}(\frac{\pi(\lvert f_{i}\rvert-f_{idl})}{2(f_{idl}-f_{min})}), & f_{min}{<}\lvert f_{i}\rvert{\leq}f_{idl}\\
\frac{k_{h_{2}}}{2}(\lvert f_{i}\rvert-f_{idl})^{2}+\frac{(\lvert f_{i}\rvert-f_{idl})^{2}}{(\lvert f_{i}\rvert-f_{max})},& \lvert f_{i}\rvert>f_{idl}
\end{cases}
\end{IEEEeqnarray*}
as it performs exactly the action described in the last sentence.
%\begin{IEEEeqnarray*}{L}
%h(f_{i})=
%\begin{cases}
%k_{h_{1}}tan^{2}(\gamma_{1}\lvert f_{i}\rvert+\gamma_{2}), & f_{min}<\lvert f_{i}\rvert\leq f_{rest}\\
%\frac{k_{h_{2}}}{2}(\lvert f_{i}\rvert-f_{rest})^{2}+\frac{(\lvert f_{i}\rvert-f_{rest})^{2}}{(\lvert f_{i}\rvert-f_{max})},& \lvert f_{i}\rvert>f_{rest}
%\end{cases}
%\end{IEEEeqnarray*}
%where $\gamma_{1}=\frac{\pi}{2(f_{idl}-f_{min})}$.%, $\gamma_{2}=-\gamma_{1}f_{idl}$.
The minimum, idle, and maximum thrusts are given by $f_{min},f_{idl},f_{max}\in\mathbb{R}^{+}$ respectively and $k_{h_{1}},k_{h_{2}}\in\mathbb{R}^{+}$.
Through the definition of $h(f_{i})$ the actuator constraints objective, implicitly has a higher priority than the position tracking objective, because $h(f_{i}){\to}\infty$ if $f_{i}{\to}{}f_{min}$ or $f_{i}{\to}{}f_{max}$.
Consequently the position tracking objective is realized only in the margins allowed by the actuator constraints.

The second component of (\ref{xi}) projects to the null-space a desired expression for the thrust magnitude $f_p$, which tracks a desired quadrotor position.
The thrust magnitude, $f_p$, has a very similar structure to (\ref{f}).
They differ only in that $f_p$ is pre-multiplied by a gain matrix and also a different gain, $k_{\xi}$, multiplies the position surface $\mathbf{s}_{x}$.
%The gain matrix assigns different weights to each axis of $f_p$ so a user can penalize the effect of the error-term of each individual axis depending on the desired maneuver, while the gain, $k_{\xi}$, is needed to adjust the influence of $\mathbf{s}_{x}$ term because as mentioned above, position tracking is performed strictly in the margins allowed by the actuator constraints.
The gain matrix assigns different weights to each axis of $f_p$ so a user can penalize the effect of the error-term of each individual axis depending on the desired maneuver, while the gain, $k_{\xi}$, is needed to adjust the influence of $\mathbf{s}_{x}$ term because position tracking is performed strictly in the margins allowed by the actuator constraints.
%The gain matrix assigns different weights to each axis of $f_p$ so a user can penalize the effect of the error-term of each individual axis depending on his/her preference while the gain, $k_{\xi}$, is needed to adjust the influence of $\mathbf{s}_{x}$ term because as mentioned above, position tracking is performed strictly in the margins allowed by the actuator constraints.
Hence it is advised that the desired position command should be in the neighborhood of the quadrotors position at the beginning of the attitude maneuver. 

Controllers to track a desired altitude command only are used in \cite{geomquadlee},\cite{geommac},\cite{geomquadlee_asian} enforcing equal priority between the attitude tracking task and the altitude tracking task.
Also no considerations for thrust saturation is mentioned, a parameter that affects both the attitude/altitude tracking tasks.
%in conjunction with (\ref{eq:distribution}) 
Here, 
%in contrast to \cite{geomquadlee},\cite{geommac},\cite{geomquadlee_asian}, 
because both the position and actuator constraint objectives are projected through $\boldsymbol{\xi}$ to the null-space of $\mathbf{A}^{\#}$ it is ensured that the attitude control objective is unobstructed assuring that the guarantees, i.e., notions of stability and regions of attraction, produced by the rigorous stability proofs hold.

%(depending on the range of the motor thrusts the desired maneuvers should not exceed a given velocity for example)
In this way, for reasonable rate maneuvers, (\ref{f_con}) always holds.
A desired maneuver is characterized as one of \textit{reasonable rate} if it is realizable in the margin of the quadrotor motor thrusts.
For example if the motor thrusts have a maximum range of 20N, a reference maneuver should not exceed a desired angular velocity in the neighborhood of $2\pi\frac{rad}{s}$.
Moreover the above solution is extremely fast to compute and implement in real time because in contrast to other works that use on-line constraint optimization algorithms (see for example \cite{Mellinger}) requiring powerful computational machinery here $\mathbf{A}^{\#}$, ${\nabla}_{\mathbf{F}}H(\mathbf{F})$, and the inverse matrix in the second component of (\ref{xi}) can be computed in an analytic form off-line.
Thus during implementation the microcontroler only needs to evaluate the precomputed analytic expressions. 
%The proposed strategy abides to controller demands (\ref{att_contr}), by using (\ref{nullspa}), because both the position and actuator constraint objectives do not interfere with the attitude control objective.

Summarizing, the strategy for avoiding thrust saturation, (\ref{nullspa}) is used during the attitude controlled mode.
Suitable desired trajectories are defined during the position controlled mode that will not require aggressive attitude deviations through polynomials of appropriate order.
Thus the thrusts remain within bounds and simultaneously position tracking is implicitly achieved as a secondary objective.
The effectiveness of the proposed solution will be verified in the next section.

\section{Simulations}
The effectiveness of the developed GNCS in performing aggressive maneuvers is verified through simulations and by a comparison to the GNCS \cite{geomquadlee}.
This GNCS was selected because it is the first quadrotor control system developed directly on SE(3) and demonstrated remarkable results in aggressive maneuvers.
%It is accepted in the literature and used in many of the most intriguing applications regarding quadrotors.
% (see \cite{Mellinger}).
%Finally \cite{geomquadlee} has a PD structure and it is easy to tune for a clear comparison with our GNCSs.
Finally \cite{geomquadlee} has similarities in structure to the proposed GNCS and is easy to tune it for a comparison.
Our focus will not only be to showcase the ability of the proposed controllers in performing aggressive maneuvers, but also to underline the ability of the developed GNCS in producing thrusts that are realizable by standard out-runner motors.

To analyze GNCSs consisting of different structure and strategies a criterion is needed for a commensurate comparison of their performance.
To this end the Root-Mean-Square (RMS) of the thrusts is used as a criterion by tuning the controllers to produce the \textit{same control effort} for a given goal.
\begin{IEEEeqnarray}{C}
f_{RMS}(t)=\sqrt{\frac{1}{t}\int_{0}^{t}\sum_{1}^{4} [f_{i}(t)]^{2} d\tau}
\label{rms}
\end{IEEEeqnarray}
Then by comparing their performance, the controller with the least error is deemed superior.

Two sets of results will be presented.
In the first set a clean comparison of the controllers without the influence of saturations or any additional factors is conducted to conclude the competence of the proposed controller.
In the second set a comparison of the GNCSs as complete solutions in the presence of saturations during an aggressive maneuver is done to evaluate the performance of the proposed GNCS.
%(it includes the strategy)Finally because the proposed GNCS includes also the strategy from Section \ref{overA}, during the attitude mode the controller performance will not be compared.
%Instead the GNCSs responses as a whole/solution with respect to the quadrotor overall performance on the maneuver will be evaluated.
The system parameters are:
\begin{IEEEeqnarray*}{C}
\mathbf{J}=[0.072,0,0;0,0.0734,0;0,0,0.1477]\;kg\,m^{2}\IEEEnonumber\\
m=1.34\;kg,
d=0.30\;m,
b_{T}=9.001\cdot10^{-3}\;m\IEEEnonumber
\end{IEEEeqnarray*}
As mentioned above the gains were tuned using (\ref{rms}) in the following manner.
First the attitude gains were tuned for a desired pitch command of $90^{o}$ followed by tuning the position gains for a desired $\mathbf{x}_{d}=[1;1;1][cm]$.
By tuning the attitude controllers first, it was ensured that during the position mode both attitude controllers embedded in their respective position control loops will produce similar/identical control efforts.
% because the desired trajectories will be designed using smooth polynomials with IC's equal to the systems state values.
%As mentioned above during the tuning procedure the saturation constraints to the thrusts were absent.
The proposed controller parameters are given by: 
\begin{IEEEeqnarray*}{C}
k_{\omega}{=}150,
k_{R}{=}5625,
\eta{=}0.809261,k_{h_{1}}{=}2,k_{h_{2}}{=}3\\
k_{v}{=}60,
k_{x}{=}900,
a{=}0.5540514,\iota_{x}{=}\iota_{y}{=}1,\iota_{z}{=}2.3,k_{\xi}{=}0.0028
\end{IEEEeqnarray*}
The benchmark controller parameters for \cite{geomquadlee} used are: 
\begin{IEEEeqnarray*}{C}
k_{\omega}=[8.6400,0,0;0,8.8080,0;0,0,17.7240]\\
k_{R}=[259.2000,0,0;0,264.2400,0;0,0,531.7200]\\
k_{v}=51.871,k_{x}=501.977
\end{IEEEeqnarray*}
%We consider two complex flight maneuvers.
%The first case corresponds to the same maneuver as in \cite{qeopidfar} for comparison
The initial conditions (IC's) are: $\mathbf{x}(0)=\mathbf{v}(0)={}^{b}\boldsymbol{\omega}(0)=0_{3\times 1},\mathbf{R}(0)=\mathbf{I}$.
The results are presented next.

\begin{figure}[!h]
\centering
%\captionsetup{font=footnotesize,singlelinecheck=false}
%\includegraphics[width=\textwidth]{./figures/Quad/frms_tun.eps}
\subfloat[\label{frmstun}]
{\includegraphics[width=0.5\columnwidth]{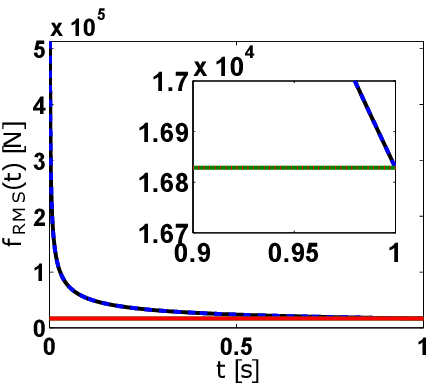}}
\put(-105,27){\parbox{\columnwidth}{${\leftarrow}1+2$}}
%\squeezeup
~
%\captionsetup{font=footnotesize,singlelinecheck=false}
\subfloat[\label{psi90tun}]
{\includegraphics[width=0.5\columnwidth]{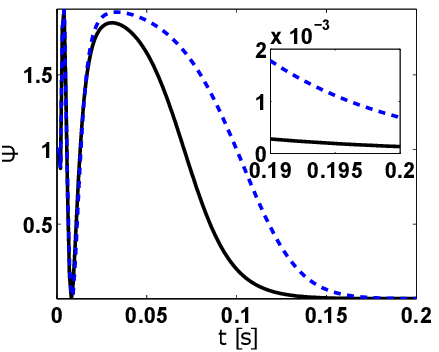}}
      \put(-79,34){\parbox{\columnwidth}{$1{\rightarrow}$}}
      \put(-41,24){\parbox{\columnwidth}{${\leftarrow}2$}}%\squeezeup
        
%\captionsetup{font=footnotesize,singlelinecheck=false}
\subfloat[\label{psitun}]
{\includegraphics[width=0.5\columnwidth]{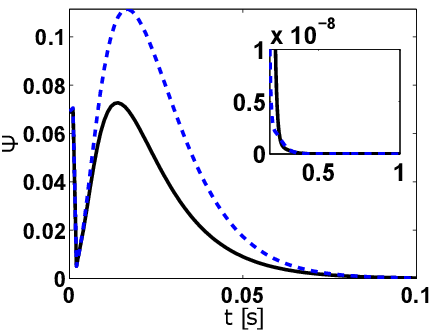}}
      \put(-88,34){\parbox{\columnwidth}{$1{\rightarrow}$}}
      \put(-75,64){\parbox{\columnwidth}{${\leftarrow}2$}}
%\squeezeup
~
%\captionsetup{font=footnotesize,singlelinecheck=false}
\subfloat[\label{trajtun}]
{\includegraphics[width=0.5\columnwidth]{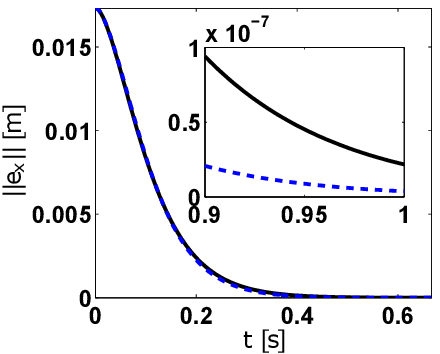}}
      \put(-70,27){\parbox{\columnwidth}{${\leftarrow}1+2$}}
%\squeezeup
        
\caption{Quadrotor response after the tuning procedure.
%(\ref{frmstun}) RMS control effort (\ref{rms}).
(\ref{frmstun}) RMS control effort calculated through (\ref{rms}).
(\ref{psi90tun}) Response for a step command of $90^{o}$.
(\ref{psitun},\ref{trajtun}) Response for a position command to $\mathbf{x}_{d}=[1;1;1][cm]$.
(\ref{psitun}) Attitude error given by (\ref{error_function}).
(\ref{trajtun})
Position error, $\lVert\mathbf{e}_{x}\rVert$.
%Position control effort calculated through (\ref{rms}).
Solid lines (1): Proposed, Dashed lines (2): Benchmark.
}
\label{Tuning}
\end{figure}

The RMS force values, (\ref{rms}), are displayed in Fig. \ref{frmstun} during the position maneuver where two overlapping horizontal lines denote the value of (\ref{rms}) at the end of the maneuver and two overlapping curved lines denote the RMS force values during the simulation.
Due to space limitations, Fig. \ref{frmstun} shows only the RMS values during the position maneuver.
%The attitude gains were tuned first and produce at the end of the attitude simulation equivalent RMS control efforts equal to 12671.159 [N], thus equal attitude control effort is exerted from both controllers.
The calculated gains produce, at the end of the attitude simulation, RMS control efforts equal to 12671.159 [N], thus equal attitude control effort is exerted from both controllers.
The RMS effort, calculated at the end of the position simulation for both controllers is $16828.54$ [N], thus equal position control effort is exerted, see Fig. \ref{frmstun}, the dashed line overlaps with the solid line and both line indexes 1, 2, point to the same curve.
The reason that the RMS control efforts are extremely large is because the control objective of precise trajectory tracking leads in the use of relatively high gains and since the controller is fed with step desired commands, extremely large control efforts are observed.
This poses no problem for an actual implementation because during trajectory tracking the control efforts remain in reasonable margins as it will demonstrated shortly.
% in the next set of results.
%Furthermore the control objective of precise trajectory tracking lead in the use of relatively high gains.
%As a result because we feed the controller with step desired commands, it is expected to observe extremely large control efforts.
%This poses no problem for an actual implementation.

Examining Fig. \ref{psi90tun}, the effectiveness of  (\ref{att_contr}) in regards to the benchmark controller is demonstrated as $\Psi$ converges to zero faster (solid black line: 1, vs dashed blue line: 2).
%In regards to the position step command Fig. (\ref{frmstun}) shows (\ref{rms}) during the position maneuver.
The quadrotor response for a position command to $\mathbf{x}_{d}{=}[1;1;1][cm]$ is shown in Fig (\ref{psitun},\ref{trajtun}).
Examining  Fig. (\ref{trajtun}), it is clear that the proposed position controller ((\ref{att_contrx}), (\ref{f})) performs equally well to the benchmark controller (solid line: 1, overlaps the dashed line: 2).
However the attitude error during the position maneuver is negotiated better by the proposed position controller as $\Psi$ converges to zero faster and with a smaller overall error, $\Psi{<}0.0727$ (solid black line: 1), vs $\Psi{<}0.1115$ (dashed line: 2), an important prevalence.%distinction.

In view of the above the ability of the proposed position controller to achieve the position command coequally to \cite{geomquadlee} while simultaneously negotiate the attitude error more efficiently makes it more effective/desirable.
%Finally the surface structure of (\ref{att_contr}), (\ref{f}) is related to our previous work \cite{VTUAV}, where in \cite{VTUAV} we proposed a controller able to negotiate bounded parametric uncertainties and disturbances.
%We reserve our investigation of this kind of properties for our proposed (\ref{att_contr}), (\ref{f}), to a future publication due to space limitations.

Having assessed the capabilities of the developed controller, a comparison between the proposed GNCS with \cite{geomquadlee} as a whole is conducted.
Thus the proposed strategy from Section \ref{overA} is now included. 
During the simulation, the following actuator constraints are in effect:
\begin{IEEEeqnarray*}{C}
f_{min}=0{ }\text{[N]}, f_{max}=20{ }\text{[N]}
\end{IEEEeqnarray*}
%\textit{(ii).} 

A complex flight maneuver is conducted, in which the quadrotor is instructed to reach a desired waypoint/position and in mid-flight perform a $360^{o}$ flip around its $\mathbf{e}_{2}$ axis.
%A complex flight maneuver will be conducted where the quadrotor is instructed to reach a desired waypoint/position and mid-flight it performs a $360^{o}$ flip around its $\mathbf{e}_{2}$ axis.
This maneuver was selected to be similar in spirit to the ones in \cite{geomquadlee},\cite{geommac},\cite{qeopidfar}, and to showcase the ability of the developed GNCS to perform aggressive maneuvers while simultaneously respecting the constraints of the out-runner motors for positive thrusts.
The initial conditions are the same as before.
%(IC's) are: $\mathbf{x}(0)=\mathbf{v}(0)={}^{b}\boldsymbol{\omega}(0)=0_{3\times 1},\mathbf{R}(0)=\mathbf{I}$.
This complex trajectory is achieved through the concatenation of the two flight modes in the following succession:
\begin{enumerate}[(a)]
\item ($t < 6$): Position Mode: At $t=0.5$ the quadrotor translates from the origin to $\mathbf{x}_{d}=[2;0;5],\mathbf{e}_{1d}=[1;0;0]$ using eighth degree smooth polynomials \cite{POLY}. 
\item ($6\leq t < 7$): Attitude Mode: The quadrotor performs a $360^{o}$ flip around its $\mathbf{e}_{2}$ axis.
$\mathbf{R}_{d}(t)$ was designed by defining the pitch angle using eighth degree polynomials.
\item ($7\leq t \leq 10$): Trajectory tracking using smooth polynomials of the same order with IC's equal to the values of the states of the quadrotor at the end of the flip and final waypoint given by $\mathbf{x}_{d}=[2;0;5],\mathbf{e}_{1d}=[1;0;0]$.
\end{enumerate}

\begin{figure}[!thpb]
\centering
%\captionsetup{font=footnotesize,singlelinecheck=false}
\subfloat[\label{PsiComp}]
{\includegraphics[width=0.5\columnwidth]{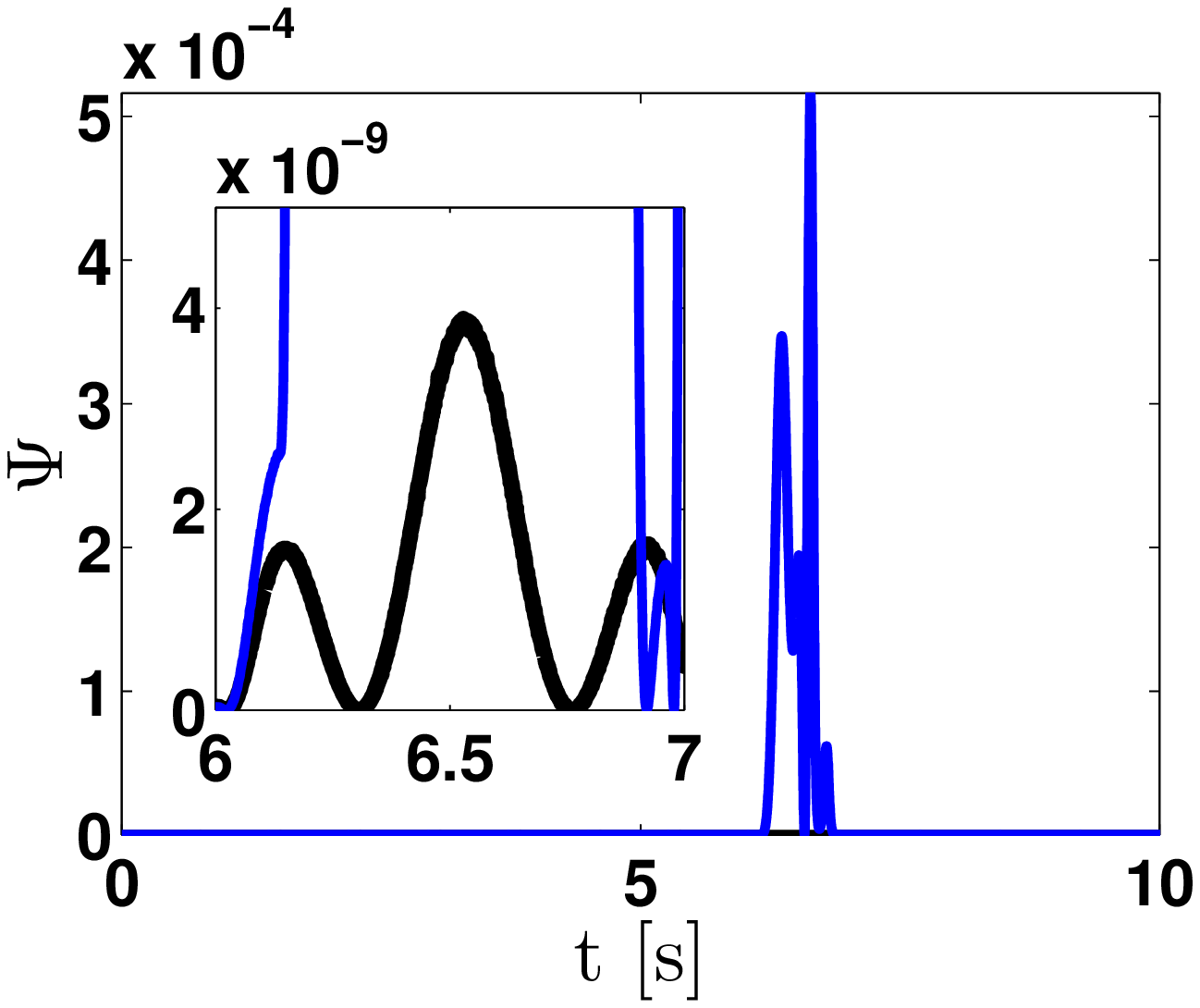}}
%\squeezeup
%\squeezeup
%\squeezeup
~
%\captionsetup{font=footnotesize,singlelinecheck=false}
\subfloat[\label{normgon}]
{\includegraphics[width=0.5\columnwidth]{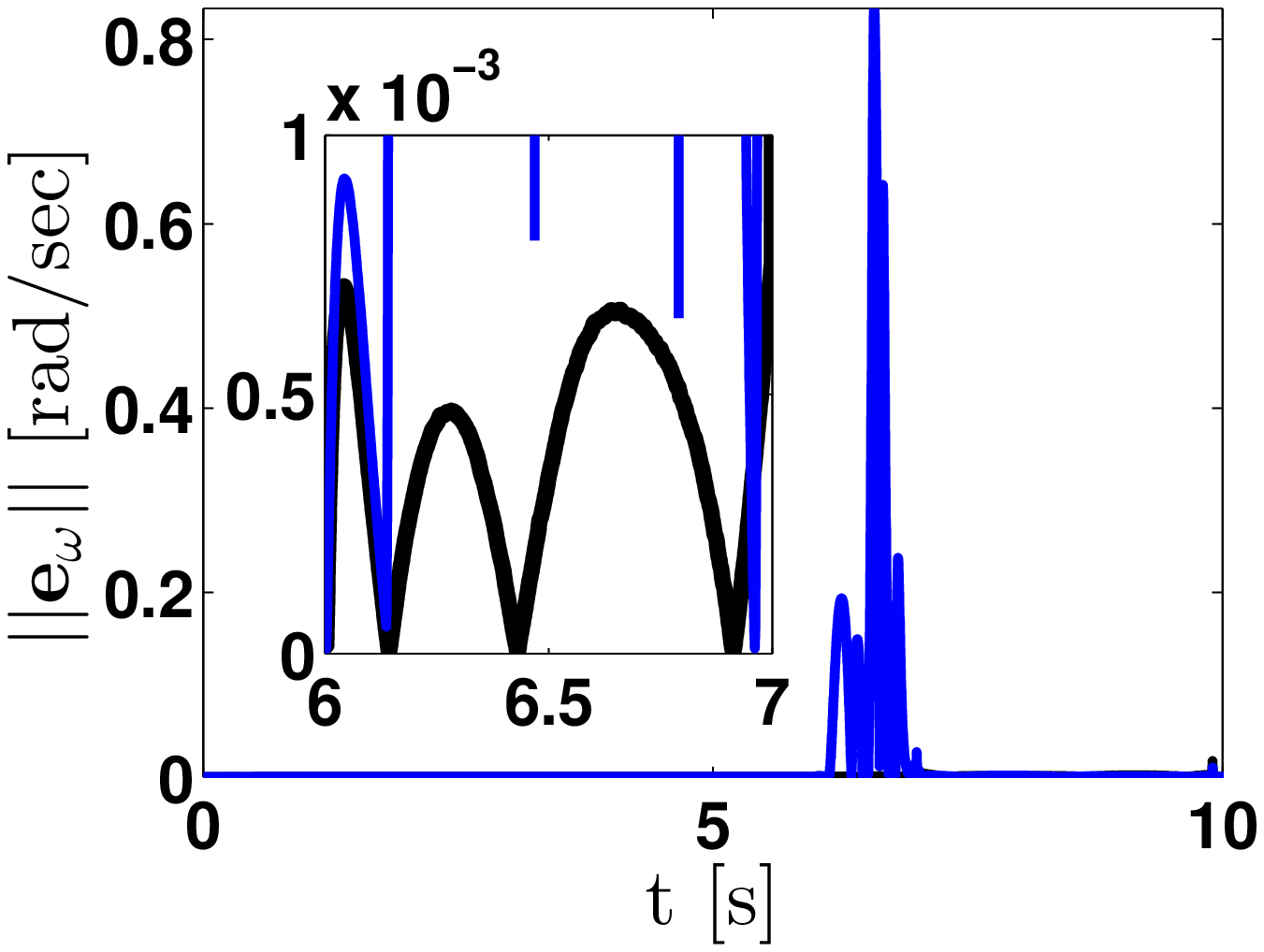}}
%\squeezeup
%\squeezeup
%\squeezeup

%\captionsetup{font=footnotesize,singlelinecheck=false}
\subfloat[\label{normpos}]
{\includegraphics[width=0.5\columnwidth]{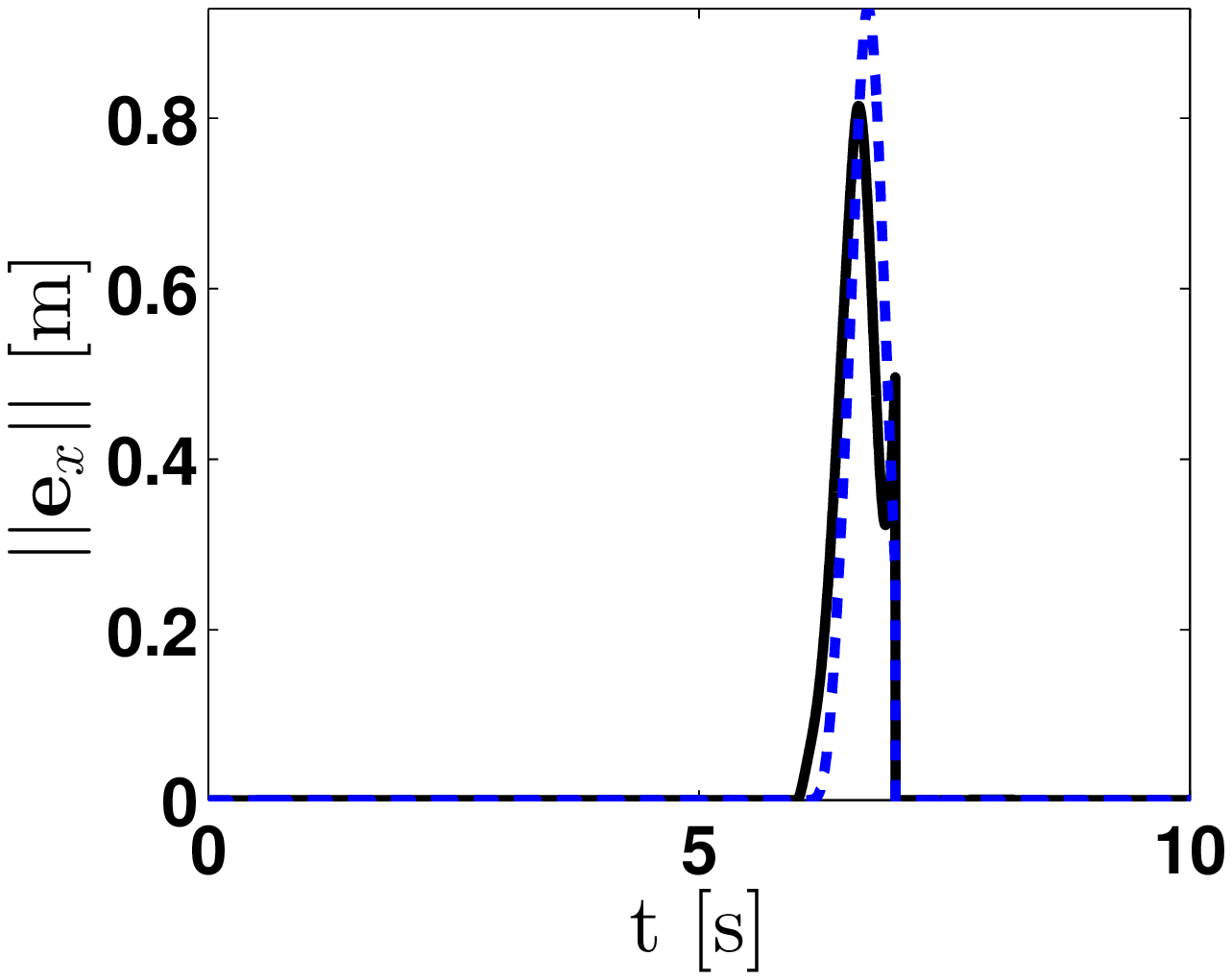}}
%\squeezeup
%\squeezeup
%\squeezeup
~
%\captionsetup{font=footnotesize,singlelinecheck=false}
\subfloat[\label{trajU}]
{\includegraphics[width=0.5\columnwidth]{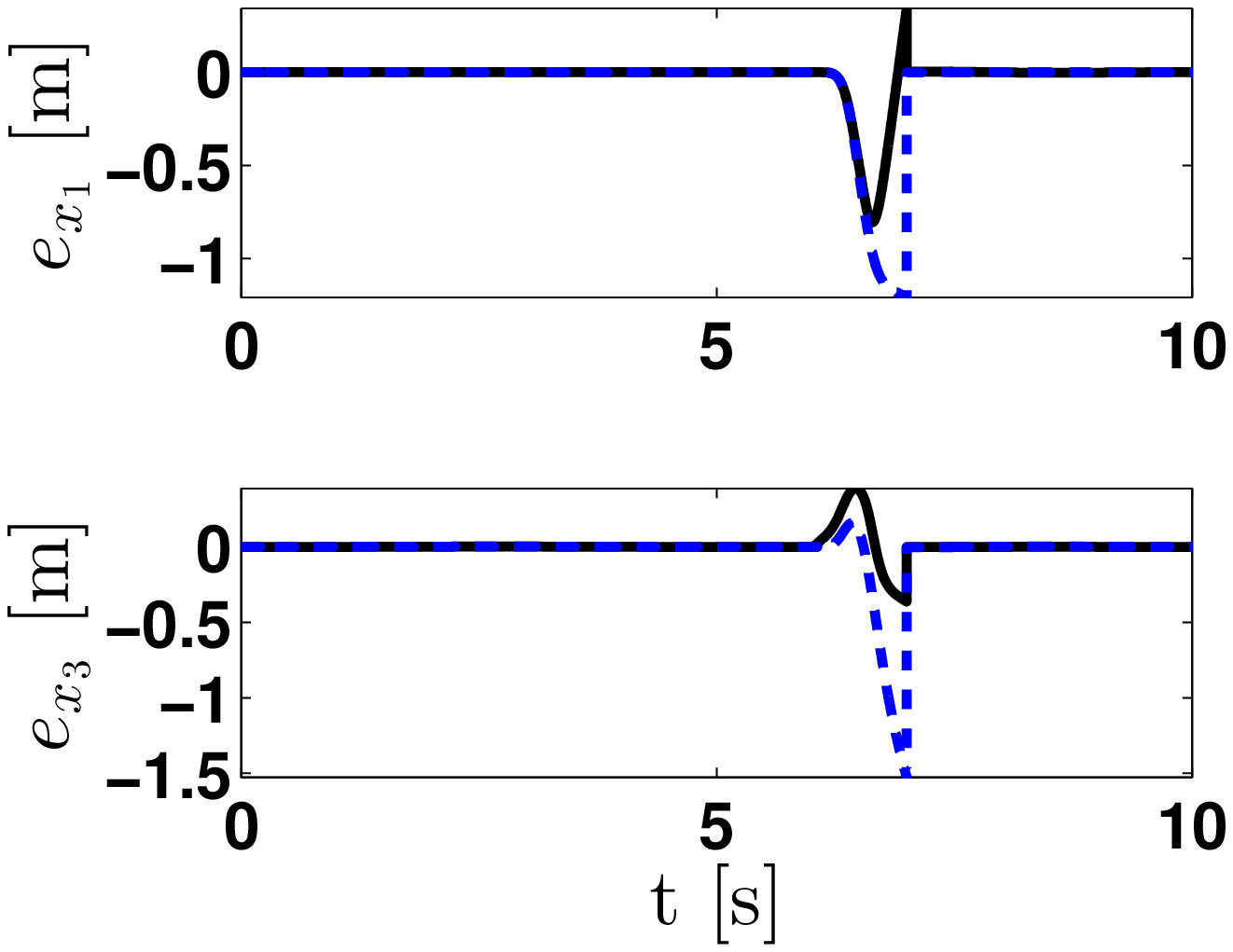}}
%\squeezeup
%\squeezeup
%\squeezeup

\subfloat[\label{thrR11}]
{\includegraphics[width=0.5\columnwidth]{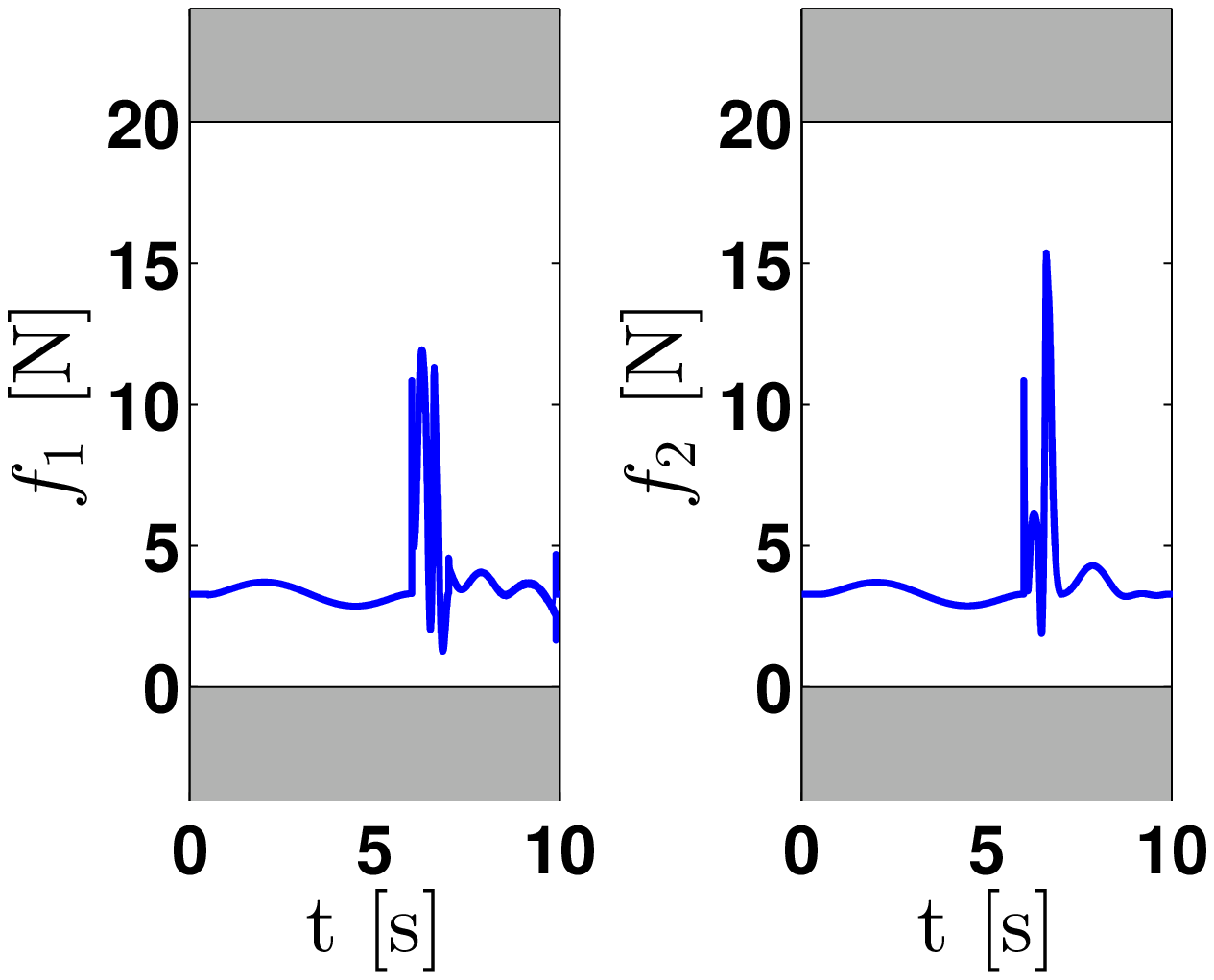}}
%\squeezeup
%\squeezeup
%\squeezeup
~
\subfloat[\label{thrR12}]
{\includegraphics[width=0.5\columnwidth]{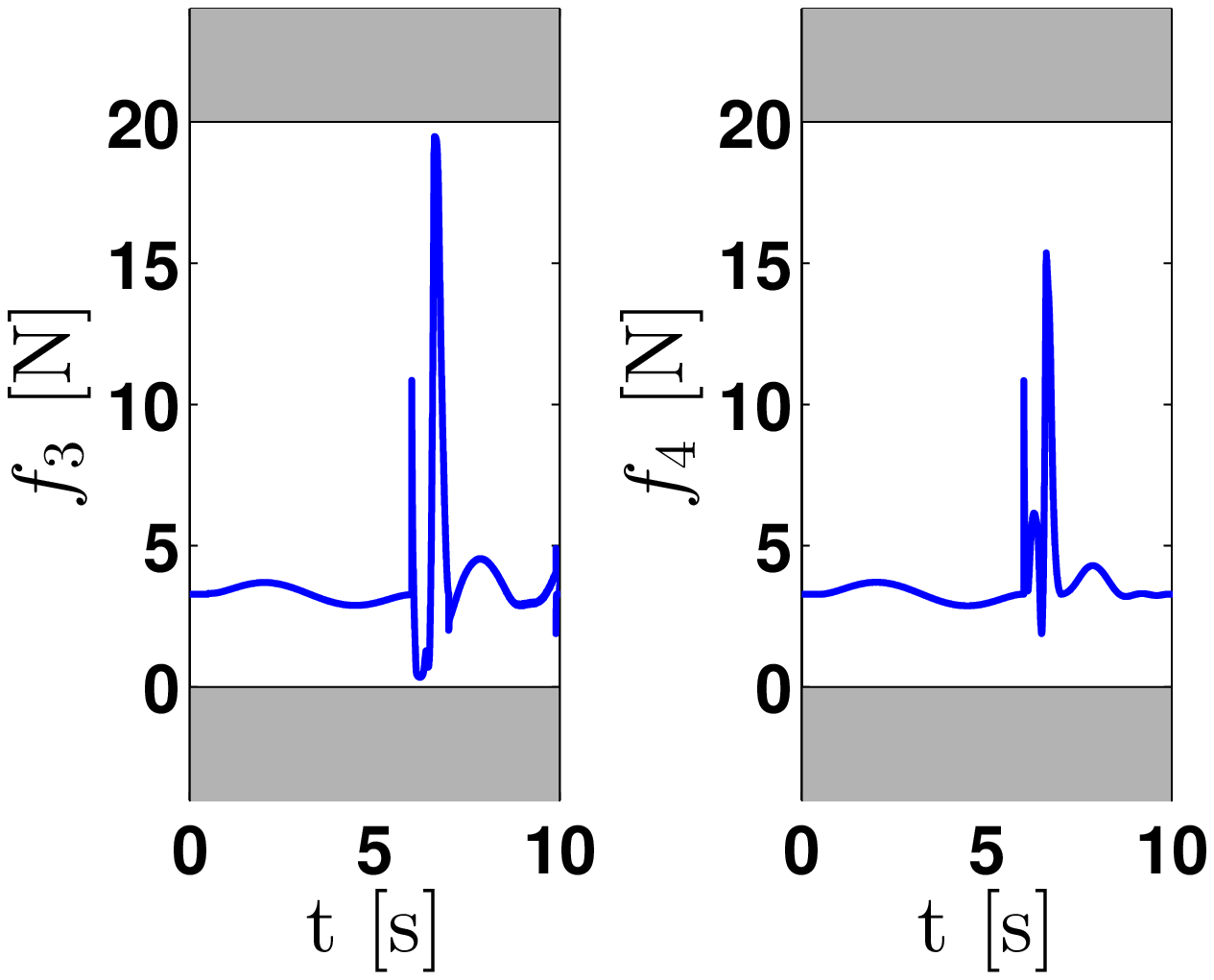}}
%\squeezeup
%\squeezeup
%\squeezeup
       
\subfloat[\label{thrM11}]
{\includegraphics[width=0.5\columnwidth]{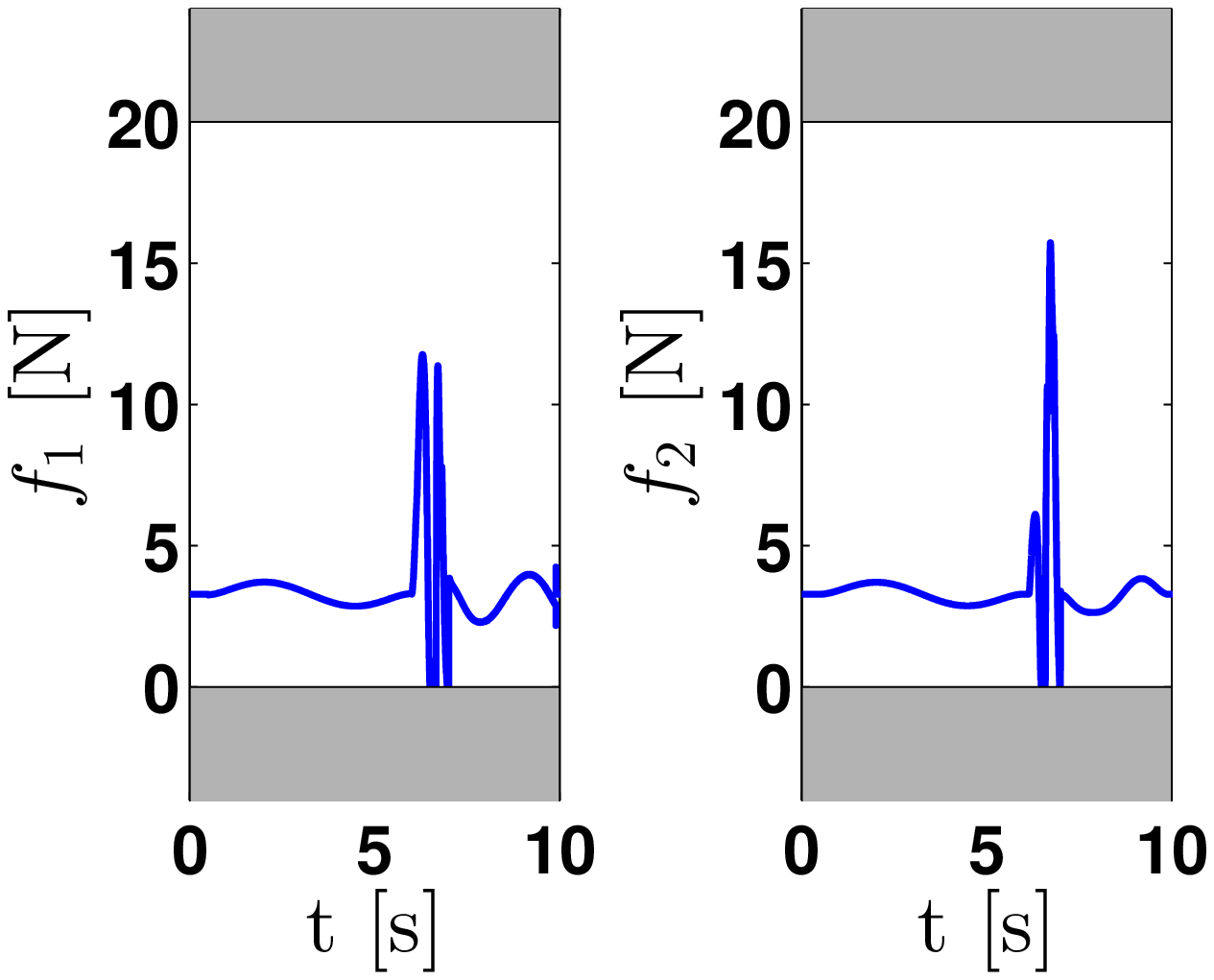}}
\put(-90,22){\parbox{\columnwidth}{${\nearrow}$}}
\put(-100,18){\parbox{\columnwidth}{sat}}
\put(-27,22){\parbox{\columnwidth}{${\nearrow}$}}
\put(-37,18){\parbox{\columnwidth}{sat}}
%\squeezeup
~
\subfloat[\label{thrM12}]
{\includegraphics[width=0.5\columnwidth]{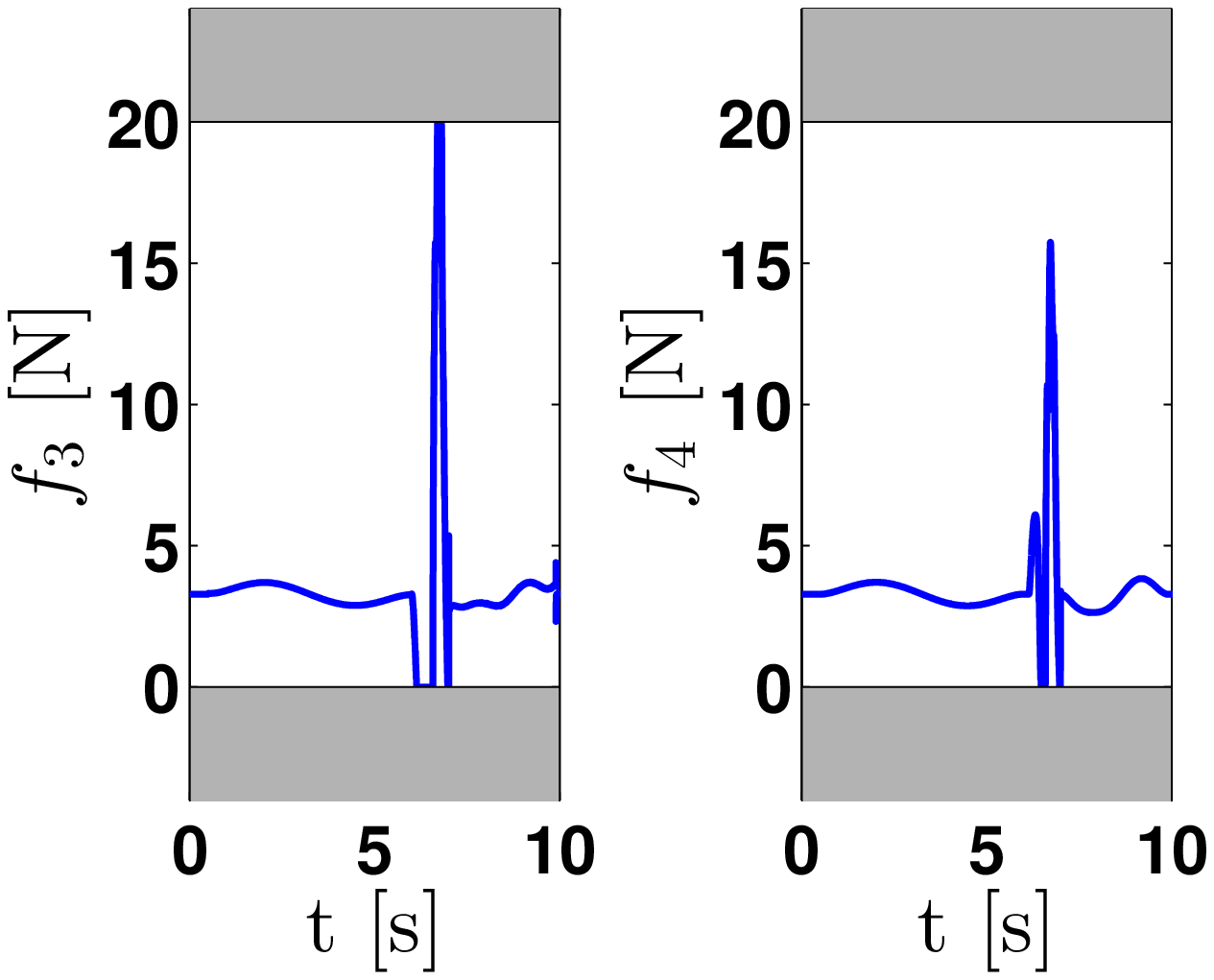}}
\put(-90,90){\parbox{\columnwidth}{${\searrow}$}}
\put(-100,94){\parbox{\columnwidth}{sat}}
\put(-90,22){\parbox{\columnwidth}{${\nearrow}$}}
\put(-100,18){\parbox{\columnwidth}{sat}}
\put(-27,22){\parbox{\columnwidth}{${\nearrow}$}}
\put(-37,18){\parbox{\columnwidth}{sat}}
%\squeezeup  

\caption{
Complex trajectory tracking.
(\ref{PsiComp}) Attitude error given by (\ref{error_function}). 
(\ref{normgon}) Angular velocity error, $\lVert\mathbf{e}_{\omega}\rVert$.
(\ref{normpos}) Position error, $\lVert\mathbf{e}_{x}\rVert$.
(\ref{trajU}) Trajectory response comparison with $f_{p}{=}0$.
(\ref{thrR11},\ref{thrR12}) Thrusts (Proposed).
(\ref{thrM11},\ref{thrM12}) Thrusts (Benchmark).}
\label{Aggressive}
\end{figure}

Simulation results of the maneuver are illustrated in Fig. \ref{Aggressive}.
%, and demonstrate the proposed controller effectiveness for precise trajectory tracking.
% with the benchmark controller to perform coequally well.
Examining Fig. \ref{PsiComp}, \ref{normgon}, \ref{normpos} it is observed that the attitude error $\Psi$, the angular velocity error $\lVert\mathbf{e}_{\omega}\rVert$, and the position error $\lVert\mathbf{e}_{x}\rVert$, only increase during the attitude portion of the maneuver (see Fig. \ref{PsiComp}, \ref{normgon}, \ref{normpos}, $6{\leq}t{<}7$).
Specifically the proposed GNCS, during the $360^{o}$ flip maneuver, demonstrates an increase only in the position tracking error, $\lVert\mathbf{e}_{x}\rVert{<}0.8142$ [m] (see Fig. \ref{normpos}, $6{\leq}t{<}7$).
The attitude error from the proposed controller remains below $\Psi{<}9.29{\cdot}10^{-9}$ ($8.36{\cdot}10^{-7}$ [deg] with respect to an axis-angle rotation) meaning that the attitude is tracked exactly, (see magnified insert in Fig. \ref{PsiComp}, thick black line), while ${}^{b}\boldsymbol{\omega}_{d}(t)$ is tracked faithfully, with $\lVert\mathbf{e}_{\omega}\rVert<0.0171 rad/s$ (see magnified insert in Fig. \ref{normgon}, thick black line).

During the same time period ($6{\leq}{t}{<}7$) the benchmark GNCS demonstrates higher tracking errors compared to the developed one.
In particular, the attitude error of the benchmark GNCS remains below $\Psi{<}5.16{\cdot}10^{-4}$ (0.046 [deg] wrt., an axis-angle rotation) denoting an error about $5{\cdot}10^{4}$ times worse compared to the developed one.
It is clear that the developed solution outperforms by far the benchmark one.
The same holds for the angular velocity error where the benchmark is competent with $\lVert\mathbf{e}_{\omega}\rVert{<}0.8337 rad/s$ (see Fig. \ref{PsiComp},\ref{normgon}, thin blue line) but exhibiting an error more than $48$ times worse.
During the $360^{o}$ flip, the benchmark position error is $\lVert\mathbf{e}_{x}\rVert{<}0.9284$ [m] (see dashed line on Fig. \ref{normpos}) signifying an error 8.37$\%$ worse.
Again the developed GNCS performs better.

The above results are attributed to the developed GNCS ability to produce thrusts that do not saturate the actuators, and still serve the attitude control objective (see Fig. \ref{PsiComp},\ref{thrR11},\ref{thrR12}).
Indeed, the lowest registered thrust is 0.3475 [N] while the largest equals to 19.4759 [N] and complies with the actuator constraints (see Fig. \ref{thrR11},\ref{thrR12}).
In contrast to this, the benchmark GNCS, during the attitude mode, is prone to thruster saturation (see `sat' in Fig. \ref{thrM11},\ref{thrM12}). 

%This signifies that our GNCSs outperforms the benchmark GNCSs and moreover our solution complies to the actuator constrains (see Fig. (\ref{thrR11},\ref{thrR12})).
During the attitude maneuver, the proposed GNCS is also able to track a desired position command, in the margins allowed by the actuator constraints, through the null-space projection of $f_{p}$ (see (\ref{xi})).
To comprehend better the effects of $f_{p}$, the same simulation was performed but this time with $f_{p}{=}0$.
The results can be seen in Fig. \ref{trajU}, where the black solid lines correspond to tracking error with active $f_{p}$, while the blue dashed lines correspond to $f_{p}{=}0$.
By comparing the responses in Fig. \ref{trajU}, it is clear that $f_{p}$ achieves its goal satisfactory well.
When $f_{p}$ is absent, the position deviation exceeds 1.2 [m] in the $x_{1}$ direction and 1.53 [m] in the $x_{3}$ direction (see Fig. (\ref{trajU}), $6{\leq}{t}{<}7$).
In contrast to this when $f_{p}$ is present, the position deviation remains below 0.8 [m] in the $x_{1}$ direction and has a mean value of close to zero in the $x_{3}$ direction (see Fig. (\ref{trajU}), $6{\leq}{t}{<}7$).
It is emphasized that the position tracking objective is achieved as a secondary task in the margins allowed by the actuator constrains.
As a result the proposed solution is intended for short durations of time.
Nevertheless the ability of the devised GNCS, to briefly track a desired position command while complying to actuator constraints, without interfering with the attitude control objective is verified.
Also it is very important to emphasize that the guarantees produced by the rigorous stability proofs hold and extend to standard quadrotors since the controller-generated thrusts produced by the GNCS are realizable through (\ref{nullspa},\ref{xi}) and comply to controller demands.
%Also it is very important to emphasize that the guarantees produced by the rigorous stability proofs hold since the controller-generated thrusts produced by the GNCS are realizable, are generated in a non invasive manner through (\ref{xi}) and comply to controller demands.
This is a very important point of this analysis, since the stability proof does not account for thrust saturations.
% and this also holds for \cite{geomquadlee},\cite{geommac},\cite{qeopidfar}.
%Thus it is very important that the desired control effort is available, i.e., avoid saturation.
%Through the simulations the proposed solution showcased results of increased precision that some might claim that precision of this magnitude is redundant, nevertheless the results are supplemented with guarantees on the system performance.
Through the simulations the proposed solution showcased results of increased precision that could be deemed redundant, nevertheless the results are supplemented with guarantees on the system performance.
% because saturations are avoided.
%Also the proposed solution (\ref{xi}), is fast and can be applied easily in real time because the pseudoinverse and all the components of (\ref{xi}) only need to be computed once and then evaluate the expresions at each iteration. % and more importantly it does not interfere with the attitude control objective.
Finally since the generated thrusts are not negative, they are realizable by standard outrunner motors and thus by the majority of quadrotors produced by the industry.
\section{Conclusion}
In this paper, new controllers for a quadrotor unmanned micro aerial vehicle were proposed based on nonlinear surfaces, composed by tracking errors that evolve directly on the nonlinear configuration manifold, inherently including in the control design the nonlinear characteristics of the SE(3) configuration space.
In particular geometric surface-based control systems were developed and through rigorous stability proofs, they were shown to have desirable closed loop properties that are almost global.
Additionally a region of attraction, independent of the position error was produced and analyzed.
A strategy allowing the quadrotor to achieve precise attitude tracking while following a desired position command and comply to actuator constraints in a computationally inexpensive manner was introduced.
This is a novel contribution differentiating this work from existing GNCSs, since the controller-generated thrusts can be realized by the majority of quadrotors available in the market.
The new features of this work were illustrated by numerical simulations of aggressive maneuvers, validating the effectiveness and the capabilities of the proposed GNCSs.

% or



%\bibliographystyle{IEEEtran}

\appendices
\section{\label{appA}}
%\subsubsection{Attitude tracking errors}
The attitude tracking errors associated with the attitude error function studied in \cite{geom_contr_book},\cite{inve_quil},\cite{qeoadapclee}, and related   properties are summarized next.

\textbf{Proposition \arabic{Prop1}.}
In regards to  (\ref{error_function}), (\ref{att_error}), (\ref{ang_vel_error}), for a given tracking command ($\mathbf{R}_{d}$, $^{b}\boldsymbol{\omega}_{d}$) and the current attitude and angular velocity ($\mathbf{R}$, $^{b}\boldsymbol{\omega}$), the following statements hold:
\begin{enumerate}[(i)]
\item $\Psi$ is locally positive-definite about $\mathbf{R}=\mathbf{R}_{d}$ and,
\begin{IEEEeqnarray}{rCl}
\lVert\mathbf{e}_{R}(\mathbf{R},\mathbf{R}_{d})\rVert^{2}&=&(2-\Psi(\mathbf{R},\mathbf{R}_{d}))\Psi(\mathbf{R},\mathbf{R}_{d})
\label{err_eq_psi}
\end{IEEEeqnarray}
\item The left-trivialized derivative of $\Psi$ is given by,
\begin{IEEEeqnarray}{rCl}
\text{T}^{*}_{I}\text{L}_{R}(\mathbf{D}_{R}\Psi(\mathbf{R},\mathbf{R}_{d}))=\mathbf{e}_{R}
\label{left_triv_der}
\end{IEEEeqnarray}
\item The critical points of $\Psi$, where $\mathbf{e}_{R}=0$, are $\{\mathbf{R}_{d}\}\cap\{\mathbf{R}_{d}\text{exp}(\pi S(\mathbf{s})),\mathbf{s}\in\text{S}^{2}\}$.
\item A lower bound of $\Psi$ is given as follows,
\begin{IEEEeqnarray}{rCl}
\frac{1}{2}\lVert\mathbf{e}_{R}(\mathbf{R},\mathbf{R}_{d})\rVert^{2}&\leq &\Psi(\mathbf{R},\mathbf{R}_{d})
\label{low_Psi}
\end{IEEEeqnarray}
%\item Let $\psi$ be a positive constant that is strictly less than 2. 6972301680
\item Let $\psi\in\mathbb{R}^{+}$.
If $\Psi(\mathbf{R},\mathbf{R}_{d})<\psi<2$, then the upper bound of $\Psi$ is given by,
\begin{IEEEeqnarray}{rCl}
\Psi(\mathbf{R},\mathbf{R}_{d})&\leq &\frac{1}{2-\psi}\lVert\mathbf{e}_{R}(\mathbf{R},\mathbf{R}_{d})\rVert^{2}
\label{upper_Psi}
\end{IEEEeqnarray}
\end{enumerate}
\textbf{Proof of Proposition \arabic{Prop1}.} See \cite{qeoadapclee}.

The associated attitude error dynamics of (\ref{error_function}), (\ref{att_error}) and (\ref{ang_vel_error}) to be used in the subsequent control design are given next.

\textbf{Proposition \arabic{Prop2}.} For the error dynamics of (\ref{error_function}), (\ref{att_error}) and (\ref{ang_vel_error}) the following hold 
\begin{IEEEeqnarray}{rCl}
\dot{\Psi}(\mathbf{R},\mathbf{R}_{d})&=&\mathbf{e}_{R}^{T}\mathbf{e}_{\omega}\\
\dot{\mathbf{e}}_{R}&=&\mathbf{E}(\mathbf{R},\mathbf{R}_{d})\mathbf{e}_{\omega}\label{der}\label{dot_Att_Error}\\
\dot{\mathbf{e}}_{\omega}&=&{}^{b}\dot{\boldsymbol{\omega}}+\mathbf{a}_{d}\IEEEnonumber\\
&=&\mathbf{J}^{-1}\left({}^{b}\mathbf{u}-{}^{b}\boldsymbol{\omega}\times\mathbf{J}{}^{b}\boldsymbol{\omega}\right)+\mathbf{a}_{d}\IEEEyesnumber\IEEEeqnarraynumspace
\label{att_error_dyn}
\end{IEEEeqnarray}
where $\mathbf{E}(\mathbf{R},\mathbf{R}_{d})\in\mathbb{R}^{3\times3}$ and $\mathbf{a}_{d}\in\mathbb{R}^{3}$ are given next,
\begin{IEEEeqnarray}{rCl}
\mathbf{E}(\mathbf{R},\mathbf{R}_{d})&=&\frac{1}{2}\{ tr [ \mathbf{R}^{T}\mathbf{R}_{d} ] \mathbf{I}-\mathbf{R}^{T}\mathbf{R}_{d} \} \\
\mathbf{a}_{d}&=&S({}^{b}{\boldsymbol{\omega}})\mathbf{R}^{T}\mathbf{R}_{d}{}^{b}{\boldsymbol{\omega}}_{d}-\mathbf{R}^{T}\mathbf{R}_{d}{}^{b}{\dot{\boldsymbol{\omega}}}_{d}
\label{E_ad}
\end{IEEEeqnarray}
with the associated norm of the matrix $\mathbf{E}(\mathbf{R},\mathbf{R}_{d})$ to be,
\begin{IEEEeqnarray}{rCl}
\lVert\mathbf{E}(\mathbf{R},\mathbf{R}_{d})\rVert&\leq&1
\label{E_norm}
\end{IEEEeqnarray}
Thus using (\ref{der}) and (\ref{E_norm}) the following inequality is satisfied,
\begin{IEEEeqnarray}{rCl}
\lVert\dot{\mathbf{e}}_{R}\rVert&\leq&\lVert\mathbf{e}_{\omega}\rVert
\end{IEEEeqnarray}
\textbf{Proof of Proposition \arabic{Prop2}.} See \cite{geomquadlee_asian},\cite{qeoadapclee}.

\section{\label{appB}}
\textbf{Proof of Proposition \arabic{Prop4}.} 
A sliding methodology is utilized through the definition of the surface in terms of the error vectors defined in (\ref{pos_error}), followed by Lyapunov analysis.
\begin{enumerate}[(a)]
\item Boundedness of $\mathbf{e}_{R}$:
The assumptions of Proposition \arabic{Prop4} imply compliance to Proposition \arabic{Prop3} thus the properties of (\ref{att_contr}) still apply.
% (intersection D L).
Equation (\ref{ep_0}) in (\ref{ater_lyap}) results in,
\begin{IEEEeqnarray}{C}
\eta k_{R} \Psi(\mathbf{R}(t),\mathbf{R}_{x}(t)){\leq} V_{\Psi}(t){\leq} V_{\Psi}(0){<}\eta k_{R}\psi_{p} \label{hkrko_Psi_p}
\end{IEEEeqnarray}
signifying that the attitude error function is bounded by,
\begin{IEEEeqnarray}{C}
\Psi(\mathbf{R}(t),\mathbf{R}_{x}(t))\leq \psi_{p} < 1,\forall t\geq 0\label{Bp_Psi}
\end{IEEEeqnarray}
\item Position Error Dynamics: The analysis that follows is developed in the following domain,
\begin{IEEEeqnarray}{rCl}
D&=&\{(\mathbf{e}_{x},\mathbf{e}_{v},\mathbf{e}_{R},\mathbf{e}_{\omega})\in\mathbb{R}^{3}\times\mathbb{R}^{3}\times\mathbb{R}^{3}\times\mathbb{R}^{3}|\IEEEnonumber\\
&{ }&\Psi(\mathbf{R},\mathbf{R}_{x})\leq\psi_{p}<1\}
\label{D}
\end{IEEEeqnarray}
\textbf{Proposition \Alph{sub}.} For initial conditions in (\ref{D}), the cosine between $\mathbf{R}\mathbf{e}_{3}$ and $\mathbf{R}_{x}\mathbf{e}_{3}$ is given by $(\mathbf{R}_{x}\mathbf{e}_{3})^{T}\mathbf{R}\mathbf{e}_{3}$ and the following holds,
\begin{IEEEeqnarray}{C}
(\mathbf{R}_{x}\mathbf{e}_{3})^{T}\mathbf{R}\mathbf{e}_{3}\geq1-\Psi(\mathbf{R},\mathbf{R}_{x})>0\label{cosine}
\end{IEEEeqnarray}
The sine of the angle between $\mathbf{R}\mathbf{e}_{3}$ and $\mathbf{R}_{x}\mathbf{e}_{3}$ is given by $( (\mathbf{R}_{x}\mathbf{e}_{3})^{T}\mathbf{R}\mathbf{e}_{3})\mathbf{R}\mathbf{e}_{3}-\mathbf{R}_{x}\mathbf{e}_{3}$ and using (\ref{err_eq_psi}),
\begin{IEEEeqnarray}{C}
\lVert((\mathbf{R}_{x}\mathbf{e}_{3})^{T}\mathbf{R}\mathbf{e}_{3})\mathbf{R}\mathbf{e}_{3}-\mathbf{R}_{x}\mathbf{e}_{3}\rVert\leq\lVert\mathbf{e}_{R}\rVert\label{b_sine}\\
\lVert\mathbf{e}_{R}\rVert=\sqrt{\Psi(2-\Psi)}\leq\sqrt{\psi_{p}(2-\psi_{p})}= \theta<1\label{b_er}
\end{IEEEeqnarray}
\textbf{Proof of Proposition \Alph{sub}.} See \cite{geomquadlee_asian},\cite{qeoadapclee}.

Equation (\ref{cosine}) is used  by adding and subtracting 
%$\frac{f\mathbf{R}_{x}\mathbf{e}_{3}}{(\mathbf{R}_{x}\mathbf{e}_{3})^{T}\mathbf{R}\mathbf{e}_{3}}$
${f\mathbf{R}_{x}\mathbf{e}_{3}}{((\mathbf{R}_{x}\mathbf{e}_{3})^{T}\mathbf{R}\mathbf{e}_{3})^{-1}}$ in (\ref{eq:position}) to obtain,
\begin{IEEEeqnarray}{rCL}
m\dot{\mathbf{v}}&=&-m\frac{k_{x}}{k_{v}}\mathbf{e}_{v}-a\mathbf{s}_{x}+\mathbf{X}+m\ddot{\mathbf{x}}_{d}\IEEEyesnumber\label{eq:refposition}
\end{IEEEeqnarray}
where $f\in\mathbb{R}$, $\mathbf{X}\in\mathbb{R}^{3}$ are given by,
\begin{IEEEeqnarray}{rCL}
f&=&\lVert \mathbf{U} \rVert (\mathbf{R}_{x}\mathbf{e}_{3})^{T}\mathbf{R}\mathbf{e}_{3}\\
\mathbf{X}&=&\lVert \mathbf{U} \rVert\left( (\mathbf{R}_{x}\mathbf{e}_{3})^{T}\mathbf{R}\mathbf{e}_{3})\mathbf{R}\mathbf{e}_{3}-\mathbf{R}_{x}\mathbf{e}_{3}\right)\\
\mathbf{U}&=&mg\mathbf{E}_{3}-m\frac{k_{x}}{k_{v}}\mathbf{e}_{v}-a\mathbf{s}_{x}+m\ddot{\mathbf{x}}_{d}\label{U}
\end{IEEEeqnarray}
Then by taking the time derivative of (\ref{pos_error}), the error dynamics of $\mathbf{e}_{v}$ are given by,
\begin{IEEEeqnarray}{rCL}
m\dot{\mathbf{e}}_{v}&=&-m\frac{k_{x}}{k_{v}}\mathbf{e}_{v}-a\mathbf{s}_{x}+\mathbf{X}\label{pos_err_dyn}
\end{IEEEeqnarray}
\item Translational dynamics Lyapunov candidate: We define,
\begin{IEEEeqnarray}{rCl}
V_{x}&=&\frac{m}{2k_{v}}\mathbf{s}_{x}^{T}\mathbf{s}_{x}+ak_{x}k_{v}\mathbf{e}^{T}_{x}\mathbf{e}_{x}
\label{pos_lyap}
\end{IEEEeqnarray}
Differentiating (\ref{pos_lyap}) and substituting (\ref{pos_err_dyn}) we get,
\begin{IEEEeqnarray}{rCl}
\dot{V}_{x}&=&\mathbf{s}^{T}_{x}(-a\mathbf{s}_{x}+\mathbf{X})+2ak_{x}k_{v}\mathbf{e}^{T}_{x}\mathbf{e}_{v}\label{Dpos_lyap}
\end{IEEEeqnarray}
Using (\ref{b_sine}-\ref{b_er}), a bound of $\mathbf{X}$ is given by,
\begin{IEEEeqnarray}{rCl}
\lVert\mathbf{X}\rVert&\leq&(B+(ak_{v}+\frac{mk_{x}}{k_{v}})\lVert\mathbf{e}_{v}\rVert+ak_{x}\lVert\mathbf{e}_{x}\rVert)\lVert\mathbf{e}_{R}\rVert\IEEEnonumber\\
&\leq&(B+(ak_{v}+\frac{mk_{x}}{k_{v}})\lVert\mathbf{e}_{v}\rVert+ak_{x}\lVert\mathbf{e}_{x}\rVert)\theta\label{xbound}
\end{IEEEeqnarray}
Defining $\mathbf{z}_{x}{=}[\lVert\mathbf{e}_{x}\rVert;\lVert\mathbf{e}_{v}\rVert]$, using  (\ref{xbound}) in (\ref{Dpos_lyap}) we arrive,
\begin{IEEEeqnarray}{rCl}
\dot{V}_{x}&\leq&-\mathbf{z}^{T}_{x}\mathbf{\Pi}_{1}\mathbf{z}_{x}+\mathbf{z}^{T}_{x}\mathbf{\Pi}_{2}\mathbf{z}_{R}\label{Dpos_Lyap}
\end{IEEEeqnarray}
and by (\ref{theta}), $\mathbf{\Pi}_{1}$ is positive definite.
\item Lyapunov candidate for the complete system: We define,
\begin{IEEEeqnarray}{rCl}
V_{g}&=&V_{x}+V
\label{glo_lyap}
\end{IEEEeqnarray}
and using (\ref{low_Psi}-\ref{upper_Psi}), (\ref{glo_lyap}) is bounded as follows,
\begin{IEEEeqnarray}{C}
\mathbf{z}^{T}_{R}\mathbf{W}_{1}\mathbf{z}_{R}{+}\mathbf{z}^{T}_{x}\mathbf{\Pi}_{3}\mathbf{z}_{x}{\leq}V_{g}{\leq} \mathbf{z}^{T}_{R}\mathbf{W}_{2}\mathbf{z}_{R}{+}\mathbf{z}^{T}_{x}\mathbf{\Pi}_{4}\mathbf{z}_{x}\label{glo_lyapb}\\
{\mathbf{\Pi}_{3}}{=}
{\begin{bmatrix}
ak_{x}k_{v}{+}\frac{mk_{x}^{2}}{2k_{v}}&{-}\frac{mk_{x}}{2}\\
{-}\frac{mk_{x}}{2}&\frac{mk_{v}}{2}
\end{bmatrix}},\;
{\mathbf{\Pi}_{4}}{=}
{\begin{bmatrix}
ak_{x}k_{v}{+}\frac{mk_{x}^{2}}{2k_{v}}&\frac{mk_{x}}{2}\\
\frac{mk_{x}}{2}&\frac{mk_{v}}{2}
\end{bmatrix}}\IEEEnonumber
\label{glo_lyapbm}
\end{IEEEeqnarray}
and both $\mathbf{\Pi}_{3},\mathbf{\Pi}_{4}$ matrices are positive definite.
Using (\ref{Datt_lyap}) and (\ref{Dpos_Lyap}) the derivative of (\ref{glo_lyap}) is given by,
\begin{IEEEeqnarray}{rCl}
\dot{V}_{g}&\leq&-\mathbf{z}^{T}_{x}\mathbf{\Pi}_{1}\mathbf{z}_{x}+\mathbf{z}^{T}_{x}\mathbf{\Pi}_{2}\mathbf{z}_{R}-\eta\mathbf{z}^{T}_{R}\mathbf{W}_{3}\mathbf{z}_{R}\label{Dglo_Lyap}
\end{IEEEeqnarray}
\item Exponential Stability: Under the conditions (\ref{theta}-\ref{w3}) of Proposition 4 all the matrices are positive definite and for $\mathbf{z}=[\lVert\mathbf{z}_{x}\rVert;\lVert\mathbf{z}_{R}\rVert]$ equation (\ref{Dglo_Lyap}) is bounded by, 
\begin{IEEEeqnarray}{rCl}
\dot{V}_{g}&\leq&-\mathbf{z}^{T}\mathbf{\Pi}_{5}\mathbf{z},\mathbf{\Pi}_{5}{=}
\begin{bmatrix}
\lambda_{min}(\mathbf{\Pi}_{1})&-\frac{1}{2}{\lVert\mathbf{\Pi}_{2}\rVert_{2}}\\
-\frac{1}{2}{\lVert\mathbf{\Pi}_{2}\rVert_{2}}&\eta\lambda_{min}(\mathbf{W}_{3})
\end{bmatrix}\label{Dglo_Lyap1}
\end{IEEEeqnarray}
Moreover (\ref{w3}) ensures that (\ref{Dglo_Lyap1}) is negative definite.
Thus the zero equilibrium of the tracking errors of the complete system dynamics is exponentially stable.
A region of attraction is given by the domain (\ref{D_x}), and (\ref{ep_0}).
%This completes the proof.
$\blacksquare$
\item Alternative regions of exponential stability:
The Lyapunov analysis above was developed in (\ref{D_x}) without restrictions on the initial position/velocity error.
This resulted to a complicated Lyapunov analysis and a reduced region of exponential stability.
Instead if we restrict our analysis to,
\begin{IEEEeqnarray}{rCl}
D_{p}&=&\{(\mathbf{e}_{x},\mathbf{e}_{v},\mathbf{e}_{R},\mathbf{e}_{\omega})\in\mathbb{R}^{3}\times\mathbb{R}^{3}\times\mathbb{R}^{3}\times\mathbb{R}^{3}|\IEEEnonumber\\
&{ }&\Psi(0){<}\psi_{p}{<}1,\lVert\mathbf{e}_{r}\rVert<\theta,\lVert \mathbf{e}_x(0)\rVert<e_{x_{max}}\}
\label{D_xmax}
\end{IEEEeqnarray}
and bound the third order error terms that arise during the analysis using $e_{x_{max}}$ then (\ref{P}), is given by
\begin{IEEEeqnarray}{C}
\mathbf{\Pi}_{1}{=}
\begin{bmatrix}
ak_{x}^{2}(1{-}\theta)&0\\
0&ak_{v}^{2}{-}\theta(mk_{x}{+}ak_{v}^{2})
\end{bmatrix}\label{p1new}\IEEEyesnumber\\
\mathbf{\Pi}_{2}=
\begin{bmatrix}
Bk_{x}&0\\
Bk_{v}+(2ak_xk_v+\frac{mk^2_x}{k_v})e_{x_{max}}&0
\end{bmatrix}\label{p2new}
\end{IEEEeqnarray}
%1,\sqrt{\psi_{p}(2-\psi_{p})},$e_{v_{max}}$ 
Alternatively a restriction on the initial velocity error results to domain,
\begin{IEEEeqnarray}{rCl}
D_{v}&=&\{(\mathbf{e}_{x},\mathbf{e}_{v},\mathbf{e}_{R},\mathbf{e}_{\omega})\in\mathbb{R}^{3}\times\mathbb{R}^{3}\times\mathbb{R}^{3}\times\mathbb{R}^{3}|\IEEEnonumber\\
&{ }&\Psi(0){<}\psi_{p}{<}1,\lVert\mathbf{e}_{r}\rVert<\theta,\lVert \mathbf{e}_v(0)\rVert<e_{v_{max}}\}
\label{D_vmax}
\end{IEEEeqnarray}
then similarly using $e_{v_{max}}$ to bound the third order error terms, $\Pi_{1}$ is given by (\ref{p1new}) and (\ref{p2new}) changes to,
\begin{IEEEeqnarray}{C}
\mathbf{\Pi}_{2}=
\begin{bmatrix}
Bk_{x}+(2ak_xk_v+\frac{mk^2_x}{k_v})e_{v_{max}}&0\\
Bk_{v}&0
\end{bmatrix}\label{p2vnew}
\end{IEEEeqnarray}
were in both cases (\ref{theta}) is given by
\begin{IEEEeqnarray}{rCl}
\theta&{<}&\min\{\frac{ak_{v}^{2}}{ak_{v}^{2}{+}mk_{x}}\}\IEEEyesnumber\label{thetaNew}\\
\IEEEnonumber
\end{IEEEeqnarray}

\squeezeup
Note that the Lyapunov analysis continues in the same manner as in Appendix \ref{appB} with (\ref{p1new}), (\ref{p2new}), (\ref{thetaNew}) (corresponding to (\ref{D_xmax})) and (\ref{p1new}), (\ref{p2vnew}), (\ref{thetaNew}) (corresponding to (\ref{D_vmax})) being utilized instead of (\ref{P}), (\ref{theta}).
It should be noted that (\ref{thetaNew}) signifies a larger basin than (\ref{theta}) but a restriction on the initial position/velocity error was introduced and this might not be desirable in some instances. 
$\blacksquare$
\end{enumerate}

\end{document}